\documentclass[a4 paper, 12pt]{article}

\usepackage{amsfonts}
\usepackage{amssymb}
\usepackage{amsmath}
\usepackage{fancyhdr}
\usepackage{eufrak}
\usepackage{latexsym}
\usepackage{amsbsy}
\usepackage[all]{xy}
\usepackage[latin1]{inputenc}
\usepackage[T1]{fontenc}        
\usepackage{amsmath}

\newenvironment{theorem}[1]
{\vskip 2mm\noindent \textbf{Theorem#1}  \it}{\vskip 2mm}

\newenvironment{mtheorem}[1]
{\vskip 2mm\noindent \textbf{Main theorem#1}  \it}{\vskip 2mm}

\newtheorem{defn}{Definition}[section]
\newtheorem{thm}[defn]{Theorem}
\newtheorem{prop}[defn]{Proposition}
\newtheorem{lem}[defn]{Lemma}

\newenvironment{corollary}[1]
{\vskip 2mm\noindent \textbf{Corollary#1}  \it}{\vskip 2mm}

\newcommand{\proof}{\vskip 2mm \noindent {\textsc{Proof: }}\rm}

\newcommand{\fin}{\hfill{\Large$\Box$}\\}
\newcommand{\finsec}{\hfill{\Large$\Box$}}

\newcommand{\al}{\alpha}

\newcommand{\Ga}{\Gamma}
\newcommand{\ga}{\gamma}
\newcommand{\pp}{\textrm{-a.e.}}
\newcommand{\si}{\sigma}

\newcommand{\om}{\omega}
\newcommand{\epsi}{\epsilon}

\newcommand{\C}{\mathbb {C}}
\newcommand{\D}{\mathbb {D}}
\newcommand{\E}{\mathbb {E}}
\newcommand{\R}{\mathbb {R}}
\newcommand{\N}{\mathbb {N}}
\newcommand{\Z}{\mathbb {Z}}

\newcommand{\Pj}{\mathbb {C}\mathbb {P}}
\newcommand{\Prr}{\mathbb {P}}

\newcommand{\cdb}{\textsf {Card }}
\newcommand{\orb}{\textsf {O}}

\newcommand{\Vol}{\textsf {vol}}
\newcommand{\Diam}{\textsf {diam }}
\newcommand{\cd}{\textsf {Card}}

\newcommand{\wert}{\, \vert \,}

\newcommand{\Per}{{\textsf {Per}  }}
\newcommand{\graphe}{{\textsf {graph}  \, }}

\newcommand{\tal}{{\tilde \alpha}}
\newcommand{\tbe}{{\tilde \beta}}

\newcommand{\Lip}{{\rm Lip \, }}
\newcommand{\Jac}{{\rm Jac \,  }}

\newcommand{\supp}{{\rm supp \, }}

\def\abs#1{\vert #1\vert}

\def\TT{{\cal T}}

\def\EE{{\cal E}}

\def\AA{{\cal A}}

\def\BB{{\cal B}}
\def\CC{{\cal C}}
\def\PP{{\cal P}}
\def\DD{{\cal D}}

\def\FF{{\cal F}}

\def\GG{{\cal G}}
\def\JJ{{\cal J}}

\def\LL{{\cal L}}
\def\VV{{\cal V}}
\def\MM{{\cal M}}
\def\GG{{\cal G}}
\def\SS{{\cal S}}

\def\QQ{{\cal Q}}

\def\WW{{\cal W}}

\def\HH{{\cal H}}
\def\com{\ar@{}[rd]|{\circlearrowleft}}

\title {Large entropy measures for endomorphisms of $\Pj^k$}

\author{Christophe Dupont}
\date{ \today }

\begin{document}

\maketitle

\begin{abstract} Let $f$ be an holomorphic endomorphism of $\mathbb{C}\mathbb{P}^k$. We construct by using coding  techniques a class of ergodic measures as limits of non-uniform probability measures on preimages of  points. We show that they have large metric entropy, close to $\log d^k$.  We establish for them strong stochastic properties and prove the positivity of their Lyapunov exponents. Since they have large entropy, those measures are supported in the support of the maximal entropy measure of $f$. They in particular provide lower bounds for the Hausdorff dimension of the Julia set. \end{abstract}

\small{\noindent \emph{Key Words}: invariant measure, metric entropy, holomorphic dynamics.

}\normalsize

%%%%%%%%%%%%%%%%%%%%%%%%%%%%%%%%%%%%
\section{Introduction}
%%%%%%%%%%%%%%%%%%%%%%%%%%%%%%%%%%%% 

This article concerns the dynamics of holomorphic endomorphisms of
$\Pj^k$. Such a mapping is
given by $k+1$ homogeneous polynomials of the same
degree $d \geq 2$ without common zero. It defines a ramified covering
of $\Pj^k$ of degree $d^k$. Its topological entropy is $h_{top}(f) =
\log d^k$, this is a consequence of \cite{G} and \cite{MP}. Remarkably, $f$ has a unique measure $\mu$ of maximal entropy
\cite{BD2}. That measure is obtained by equidistributing preimages of generic points: there exists an algebraic subset $\EE \subset \Pj^k$ depending on $f$ such that for every $z \notin \EE$:
\[ {1 \over d^{kn}} \sum_{f^n(y) = z}  \,\delta_y \to \mu ,  \]
where $\delta_y$ denotes
the Dirac mass at $y \in \Pj^k$ (see \cite{BD2},
\cite{DS1}). Strong stochastic properties have been established for $\mu$: exponential decay of correlations,  Central Limit Theorem,  Almost Sure
Invariance Principle and Large Deviations Theorem, at least for
H\"older observables (see \cite{DNS}, \cite{D1} and references
therein). Moreover, we know that the Lyapunov exponents of $\mu$ are
larger than or equal to ${1 \over 2} \log d$ \cite{BD1}. We refer to the
survey article \cite{DS2} for more details concerning those properties.  \\

In the present article, we use coding techniques to construct for endomorphisms of $\Pj^k$ a class of ergodic measures with entropy close to $\log d^k$. Those measures are limits of  
non-uniform probability measures on $f^{-n}(z)$, and, as we shall see, their support is included in $\JJ := \supp \mu$. We prove for them the Almost Sure Invariance Principle (ASIP) for H\"older observables  and the positivity of the Lyapunov exponents. 

Let us note that coding techniques were introduced in \cite{PUZ}  for rational fractions acting on $\Pj^1$. That allowed them to construct a class of ergodic measures containing $\mu$ and satisfying the ASIP for H\"older observables. That coding method was recently extended to $\Pj^k$ in order to study the maximal entropy measure $\mu$ \cite{D1}. Here we continue to develop the coding techniques in higher dimensions. \\

 Let us introduce notations related to coding. Let $\AA :=
\{ 1, \ldots, d^k\}$, $\Sigma := \AA^{\N}$ and $s$ be the left
shift acting on $\Sigma$. We set $z$ outside the critical values of every $f^n$. A map $z_n$ will  denote a bijection $\AA^n \to f^{-n}(z)$,  extended to $\Sigma$ as a constant map on
$n$-cylinders. We say that a sequence $(z_n)_n$ is \emph{compatible} if $f \circ z_{n+1} = z_n \circ s$. Given an H\"older function $\varphi :
\Sigma \to \R$, let  $\nu_\varphi$ denote its associated Gibbs
measure on $\Sigma$. This is the unique $s$-invariant measure
satisfying $P_\varphi =  \int \varphi \, d\nu_\varphi +
h_{\nu_\varphi}(s)$, where $P_\varphi$ and $h_{\nu_\varphi}(s)$ respectively denote the pressure of $\varphi$ and the entropy of $\nu_\varphi$ (see \cite{B}, \cite{W}). \\

We shall actually work with potentials $\varphi : \Sigma \to \R$ satisfying $\tau_{\theta}(\varphi) :=
P_\varphi - \abs {\varphi}_\infty - \log d^{k-\theta} > 0$ for small $\theta$'s, less than some $\theta_k$ specified in section \ref{state}. The  measure $\nu_\varphi$ can therefore be a product measure with weights close to $d^{-k}$ (see corollary 1). Observe that  $\tau_{\theta}(\varphi) > 0$ combined with $h_{\nu_\varphi}(s) = P_\varphi - \int \varphi
\, d\nu_\varphi$ implies $h_{\nu_\varphi}(s) > \log
d^{k-\theta}$. \\

Let us state our main theorem. We set $\nu_n := (z_n)_*\nu_\varphi$, that pushforward measure is supported on $f^{-n}(z)$ and satisfy $\nu_n (y) = \sum_{f(x) = y} \nu_{n+1}(x)$ by the compatible condition and the $s$-invariance of $\nu_\varphi$. 

\begin{mtheorem}{:}
For every $\theta < \theta_k$, there exists a zero volume subset $\FF_\theta \subset \Pj^k$ satisfying the
following property. For every $z \notin \FF_\theta$, there exist compatible
bijections $z_n : \AA^n \to f^{-n}(z)$ such that, for every H\"older
function $\varphi : \Sigma \to \R$ satisfying $\tau_\theta(\varphi) >
0$:

\begin{enumerate}

\item[{\bf (a)}] $\sum_{f^n(y)=z}
  \nu_n(y) \, \delta_y$ converges to a mixing $f$-invariant measure $\nu$ on $\Pj^k$.

\item[{\bf (b)}] The metric entropy is preserved, i.e. $h_\nu(f) = h_{\nu_\varphi}(s)$.

\item[{\bf (c)}] $\nu$ does not charge the
  algebraic subsets of $\Pj^k$ and $\supp \nu \subset \JJ = \supp \mu$.

\item[{\bf (d)}] The Lyapunov exponents of $\nu$ are larger than or equal to ${1 \over 2}(h_\nu(f) - \log d^{k-1}) > 0$.  

\item[{\bf (e)}] The measure $\nu$ satisfies the exponential decay of correlations and the Almost Sure Invariance Principle for
   H\"older observables.
  
\end{enumerate}
\end{mtheorem}

 The proof relies on independent results which are stated in
section \ref{state}.   Let us make some comments. Item {\bf (b)} and condition $\tau_\theta(\varphi) > 0$ imply that $h_{\nu}(f) > \log
d^{k-\theta}$. That large entropy property is crucial for the proofs of {\bf(c)} and {\bf(d)}. We shall see that those assertions are actually more general, they hold for every ergodic measure with entropy larger than $\log d^{k-1}$. We observe also that {\bf (b)} implies $h_\nu(f) < \log d^k$ once $\nu_\varphi$ is not the maximal entropy measure on $\Sigma$. This proves that $\nu$, whose support is included in $\JJ$ by {\bf (d)}, is singular with respect to $\mu$ (these are different ergodic measures). Concerning {\bf (e)}, we recall that an observable $\psi : \Pj^k \to \R$ satisfies the ASIP if the process defined by the Birkhoff sums of $\psi$ follows $\nu$-almost everywhere  the trajectory of a Brownian motion. That strong stochastic property
implies the Central Limit Theorem and the
  Law of Iterated Logarithm for $\nu$.\\

The  following corollary concerns the case when $\nu_\varphi$ is a product measure on $\Sigma$. We set for the first item $y =
z_n(\al_0,\ldots,\al_{n-1}) \in \Pj^k$.

\begin{corollary}{ 1:} Let $\eta < \theta_k$ and  $(w_\al)_{\al \in \AA}$ satisfying $\sum w_\al =1$ and $w_\al \leq
d^{-k+\eta}$. For every $z$ outside a zero volume subset, there exist
compatible bijections $z_n$ such that:  
\begin{enumerate}
\item[-] $\sum_{f^n(y) = z} \Big(  \prod_{i=0}^{n-1} w_{\al_i} \Big) \delta_y $ converges to a mixing $f$-invariant measure $\nu$.

\item[-] The metric entropy of $\nu$ is equal to $h_\nu(f) =  - \sum w_\al \log w_\al$.

\item[-] The measure $\nu$ satisfies {\bf (c)}, {\bf (d)}  and {\bf (e)}.
\end{enumerate}
\end{corollary}
That corollary  follows by taking  $\varphi(\al_0,\al_1, \ldots) := \log w_{\al_0}$ in the main theorem. We have in that case $P_\varphi = 0$. Observe that the condition $\tau_{\theta}(\varphi) =  - \abs {\varphi}_\infty - \log d^{k-\theta} > 0$ is  then satisfied for every $\theta \in ]\eta  , \theta_k[$. \\

We deduce the following property, where $\MM$ is the set of $f$-invariant measures.

\begin{corollary}{ 2:}
The image of the entropy function $h : \MM \to [0, \log d^k]$ contains a non trivial interval of the form $[\log d^k - \epsi , \log d^k]$.
\end{corollary}
That follows from the continuity of  $-\sum w_\al \log w_\al$ and the fact that the maximal value occurs exactly when $w_\al = d^{-k}$ for every $\al \in \AA$. We observe that, since $f$ is a $C^\infty$ mapping, the entropy function $h$ is upper semicontinuous  \cite{N2}. That property ensures the existence of a maximal entropy measures for $f$, but does not provide corollary 2.  \\

The first part of item {\bf (c)} is obviously related to the size of $\nu$. Concerning Hausdorff dimension, the article \cite{D2} yields for every ergodic measure with positive exponents a lower bound for the dimension of positive $\nu$-measure Borel sets. We therefore get:
\begin{corollary}{ 3:}
Let $\nu$ be an ergodic measure provided by the main theorem. Let $\lambda_1 (\nu) \geq \ldots \geq
\lambda_k (\nu)$ denote its Lyapunov exponents (they are positive by {\bf (d)}) and $A$ be a Borel subset of positive $\nu$-measure. Then the Hausdorff dimension of  $A$ satisfies 
\begin{equation*}\label{mmj}
 \dim_\HH A \geq {\log d^{k-1} \over \lambda_1 (\nu)} + { h_\nu(f) -  \log d^{k-1} \over \lambda_k (\nu)}.
 \end{equation*}
\end{corollary}
That estimate in particular holds for $A = \JJ$, hence the ergodic measures provided by the main theorem yield lower bounds for the Hausdorff dimension of $\JJ$. Let us notice that when $k=1$,  the formula $\inf \, \{ \dim_\HH A \, , \, \nu(A) > 0 \}  = {h_\nu(f)  \over \lambda (\nu)}$ holds for every ergodic measure $\nu$ with positive entropy \cite{M}. The proof relies on the fact that rational fractions are conformal mappings. An analogous formula remains unknown in  higher dimensions. We refer to the articles \cite{BDM}, \cite{DD} and \cite{D2} for results concerning that problem. The article \cite{D2} in particular yields for the maximal entropy measure and for $k=2$ the bound $\inf \,  \{ \dim_\HH A \, , \, \mu(A) > 0 \}  \geq  {\log d \over \lambda_1(\mu) } + {  \log d \over \lambda_2(\mu)}$, which is half of the formula conjectured in \cite{BDM}. \\

Observe also that  {\bf (d)} ensures that the function $\log \Jac f$,
which is $-\infty$ on
the critical set of $f$, belongs to $L^1(\nu)$. That follows from the classical  formula 
\[\int_{\Pj^k} \log \Jac f \, d\nu= 2(\lambda_1(\nu) + \ldots + \lambda_k(\nu)) . \] 
We notice that the integrability of $\log \Jac f$ for the equilibrium measure $\mu$ can be shown without using the exponents: this is proved in \cite{FS} and \cite{DS1} by  using respectively the construction of $\mu$ from pluripotential theory and from pullbacks of volume forms.  \\

We finally note that the existence and unicity of equilibrium states was recently established for endomorphisms of $\Pj^k$ for a class of H\"older functions \cite{UZ}. The approach consists there in studying the Ruelle-Perron-Frobenius
operator. It should be interesting to compare those equilibrium states with the measures constructed here. 

%%%%%%%%%%%%%%%%%%%%%%
\section{Structure of the proof of the main theorem}\label{state}
%%%%%%%%%%%%%%%%%%%%%%

The main theorem follows from theorems A, B and C stated below. We denote by $d(\cdot,\cdot)$ the standard distance on $\Pj^k$ and we set $\theta_k := {2 \over 5(k-1)}$.

\begin{theorem}{ A:}
For every $\theta < \theta_k$, there exists a zero volume subset $\SS_\theta \subset \Pj^k$ satisfying the
following property. For every $z \notin \SS_\theta$, there exist compatible bijections $z_n : \AA^n \to f^{-n}(z)$,
 an increasing sequence of subsets $(\GG(n))_n \subset \Sigma$ and $\rho > 0$ such that, for every H\"older function $\varphi : \Sigma \to \R$ satisfying
  $\tau = \tau_\theta(\varphi) > 0$:
\begin{enumerate}
\item  $(z_n)_n$ converge
  $\nu_\varphi \pp$ to some limit map $\om : \Sigma \to \Pj^k$.

\item $\nu_\varphi(\GG(n)) \geq 1 -  c_\tau  \, e^{- n \tau}$ and $d( z_n(\tal), 
\om(\tal) ) \leq  c_\rho \, d^{-\rho n}$ for every $\tal \in \GG(n)$. 

\item The pushforward measure $\om_* \nu_\varphi$ is a mixing $f$-invariant measure satisfying the Almost Sure
  Invariance Principle and the exponential decay of correlations for
  H\"older observables.
\end{enumerate}
\end{theorem}

The point 1 implies that $(z_n)_* \nu_\varphi$, which by definition is equal to $\sum_{f^n(y)=z}
  \nu_n(y)  \, \delta_y$, converges to $\om_* \nu_\varphi$. It also implies the relation $f \circ
  \om = \om \circ s$ by using the compatible
  condition on $(z_n)_n$. In particular, since $\nu_\varphi$ is a mixing invariant measure, that property holds for $\om_* \nu_\varphi$. That proves {\bf(a)} of the main theorem for $\nu := \om_*
  \nu_\varphi$. The point 2 of theorem A is crucial to show its point 3, which
  corresponds to {\bf (e)}.

Theorem A was proved by Przytycki-Urba\'nski-Zdunik for $k=1$ and every H\"older
  continuous function $\varphi$ (see \cite{PUZ}, section 3), the proof relies there on Koebe distortion theorem. We use in the context $k \geq 2$ a quantitative version of the inverse branch Briend-Duval's theorem established in \cite{D1}: that result allows to control the size of at least $(1 - d^{-n\theta})d^{kn}$ inverse branches of $f^n$. The value $\theta_k = {2 \over 5(k-1)}$ is due to the technique of the proof, which  requires to work on complex lines and to estimate moduli of annuli.  \\
     
The following theorem  yields {\bf(b)} by setting $\FF_\theta := \SS_\theta \cup \HH$. 
\begin{theorem}{ B:}
There exists a zero volume subset $\HH \subset \Pj^k$ satisfying the
  following property. Let $\theta < \theta_k$ and $z \notin \SS_\theta \cup \HH$. Let $\varphi : \Sigma \to \R$ be an H\"older function satisfying
  $\tau_\theta(\varphi) > 0$ and let $\om : \Sigma \to \Pj^k$  provided by theorem A. Then $h_{\om_* \nu_\varphi} (f) = h_{\nu_\varphi}(s)$.
\end{theorem}

The proof consists in estimating the size of the fibers of  $\om$ and applying Abramov-Rohlin's formula for the entropy of a skew product. A similar result was proved by Przytycki (\cite {P}, section 2) in dimension $k=1$ for maps $\om$ which are boundary extensions of Riemann mappings on the unit disc. Here we follow the same approach by showing that the fibers of $\om$ are small for good choices of the root $z$, namely when the fibers $f^{-n}(z)$ are far enough from the critical set of $f$. The set $\HH$ is introduced in order to ensure that property. \\

Theorem B, combined with $\tau_\theta(\varphi) > 0$, yields $h_\nu(f) = h_{\nu_\varphi}(s) > \log d^{k-\theta}$. Let us see
how that implies  {\bf(c)}. We use for that purpose the relativized variational
principle: if $\nu$ is ergodic and if $\nu(A)
> 0$, then  $h_\nu(f)$ is less than or equal to $h_{top}(f,A)$, the topological entropy of $f$ relative to $A$ (see \cite{BD2}, section 4). We
deduce that $\nu$ does not charge the algebraic subsets of $\Pj^k$. We indeed have $h_{top}(f,A) \leq \log d^p$ for every $p$-dimensional algebraic
subset $A$, that follows from Gromov's argument \cite{G}.  Concerning
the inclusion $\supp \nu \subset \JJ$, Dinh \cite{D} proved that $h_{top}(f,A) \leq \log d^{k-1}$ for every compact set $A$ not intersecting $\JJ$. That completes the proof of {\bf(c)}. \\

The item {\bf(d)} is provided by the following general result:
\begin{theorem}{ C:}
Let $f$ be an holomorphic endomorphism of $\Pj^k$ of degree $d \geq 2$ and let $m$ be an ergodic
$f$-invariant measure. If $h_m(f) > \log d^{k-1}$, then the Lyapunov exponents of $m$ are larger than or equal to ${1 \over 2} (h_m(f) - \log d^{k-1} ) > 0$.
\end{theorem}
The proof is a consequence of theorem D below (this is explained in section \ref{cetd}). Let us denote by $\Lambda_1 > \ldots > \Lambda_q \geq
- \infty$ the exponents of $m$ and by $m_1 , \ldots , m_q$ their
multiplicities. Recall that Margulis-Ruelle inequality \cite{R} states $h_m(f) \leq  2
m_1 \Lambda_1^+ + \ldots + 2 m_q \Lambda_q^+$, where $x^+ := \max \{  x, 0\}$. The next result  extends that inequality.

\begin{theorem}{ D:}  Assume that $q \geq 2$. Then for every $2 \leq j \leq q$:
\[  h_m(f) \leq \log d^{m_1 + \ldots + m_{j-1}}  + 2 m_j \Lambda_j^+ +
\ldots + 2 m_q \Lambda_q^+ .\]
\end{theorem}

Theorem D and its corollary theorem C were established by de Th\'elin \cite{dT} assuming
$\Lambda_q >- \infty$, namely $\log \Jac f \in L^1(m)$. Our aim is to
extend it to $\Lambda_q = - \infty$. The proof of theorem D  in
\cite{dT} is based on the construction of \emph{relative} almost
stable manifolds and on  volume estimates. De Th\'elin told us that it
was possible to extend theorem D to the non-integrable case by using
the same method, here we verify it. 

The stable manifolds were obtained in \cite{dT} by composing  forward graph transforms for $f^{-1}$ along $m$-generic orbits. In the non-integrable case, we obtain them by performing instead backward graph transforms for $f$ itself: we use for that purpose the method of Hirsch-Pugh-Shub \cite{HPS} which allows  non-injective maps. The Oseledec-Pesin's charts are provided in our context by Newhouse
theorem \cite{N1}. Then, once the stable manifolds are constructed,
volume estimates are obtained by slicing
arguments as in \cite{dT}. It turns out that an occurence of multiplicities, due to a lack of injectivity, does not affect the bounds. 

We thank de Th\'elin for several discussions concerning that problem of extension. Concerning this topic, we notice that  Buzzi \cite{B} announced us an estimate similar to theorem D in a
real setting, with the term $\log d^{m_1 + \ldots + m_{j-1}}$ being replaced by the topological entropy of embedded smooth discs of dimension $m_1 + \ldots + m_{j-1}$. \\

The article is organized as follows. Classical facts concerning Gibbs measure  are recalled in section 3. Theorems A, B and C-D are proved in sections 4, 5 and 6. In the sequel $c$
denotes a constant which may change from a line to another.

%%%%%%%%%%%%%%%%%%%%%%%%%%%%%%%%%%%%%%%%%%%%%%%%%%%%%%%%%
\section{Gibbs measures on  $(\Sigma,s)$}\label{gibbs} 
%%%%%%%%%%%%%%%%%%%%%%%%%%%%%%%%%%%%%%%%%%%%%%%%%%%%%%%%%

Our references are the classical books \cite{B}, \cite{W}.  Let $\tal = (\al_i)_{i\geq
0}$ be the elements of $\Sigma = \AA^\N$. We denote by $\AA_n$ the $n$-cylinders of $\Sigma$, and by $[\al_0,\ldots,\al_{n-1}]$ a $n$-cylinder. Let also $\AA_{i,j} :=
s^{-i} \AA_{j-i}$ for $i \leq j$, observe that $\AA_n = \AA_{0,n}$. We endow $\Sigma$ with a product metric. Let $\varphi$ be an
H\"older function $\Sigma \to \R$ and $\nu_\varphi$ be its associated equilibrium measure. This is the unique $s$-invariant measure on $\Sigma$ satisfying 
$P_\varphi = \int_\Sigma \varphi \, d\nu_\varphi +
h_{\nu_\varphi}(s)$, where $P_\varphi$ is the pressure of $\varphi$ and $h_{\nu_\varphi}(s)$ is the metric entropy of $\nu_\varphi$. This is also the unique $s$-invariant measure for
which there exist $c_1, c_2 > 0$ satisfying for every $\tal \in \Sigma$ and $n \geq 1$:
\begin{equation}\label{ggiibb}
  c_1 \,
  e^{\sum_{i=0}^{n-1} \varphi \circ s^i(\tal) - n P_\varphi}   \leq \nu_\varphi [\al_0,\ldots,\al_{n-1}]  \leq c_2 \,
  e^{\sum_{i=0}^{n-1} \varphi \circ s^i(\tal) - n P_\varphi} . 
  \end{equation}
It turns out that $\nu_\varphi$ satisfies the following exponential mixing property (see \cite{B}, Proposition 1.14): there exists $c, \delta > 0$ such that for every $i \leq j < k \leq l$,
\begin{equation}\label{ibbs}
  \forall E \in \AA_{i,j} \, , \, \forall F \in \AA_{k,l} \ , \ \abs{ \nu_\varphi(E \cap F) - \nu_\varphi(E) \, \nu_\varphi(F) }  \leq c \, \nu_\varphi(E) \,  \nu_\varphi(F) \,  e^{-\delta(k-j)}.   
 \end{equation}
That property implies that $\nu_\varphi$ is mixing in the usual sense. We shall need (\ref{ibbs}) for theorems \ref{exdecay} and \ref{philstout} below. We say that an observable  $\chi : \Sigma \to \R$ is $\nu_\varphi$-centered if $\int_\Sigma \chi \, d\nu_\varphi = 0$. We denote by $\E(\chi \vert \AA_n)$ the conditional
expectation of $\chi \in L^1(\nu_\varphi)$ with respect to the partition $\AA_n$. Let us introduce the following definition.

\begin{defn}\label{ddee}
We say that an observable $\chi : \Sigma \to \R$ is
\emph{$L^p(\nu_\varphi)$-cylinder} ($p \geq 1$) if 
\[ \chi \in L^p(\nu_\varphi) \ \ \textrm{ and } \ \ \abs{  \chi -
  \E(\chi \vert  \AA_n)  }_p \leq c \, e^{-\ga n} \textrm{ for some }
\gamma > 0. \]
\end{defn}

\noindent The next result states the exponential decay of correlations for those observables. 

\begin{thm}\label{exdecay}
Let $\chi_1,\chi_2 : \Sigma \to \R$ be bounded $\nu_\varphi$-centered observables which are $L^1(\nu_\varphi)$-cylinder. Then $\left \vert \int_{\Sigma} \chi_1 \cdot
\chi_2 \circ s^n \, d\nu_\varphi \right \vert \leq  c \,
e^{-n\lambda}$ for some $\lambda > 0$.
\end{thm}

\proof Let $\chi_{j,m} := \E(\chi_j \vert
\AA_m)$ and write:
\[  \chi_1 \cdot \chi_2 \circ s^n = (\chi_1 - \chi_{1,m}) \cdot \chi_2
\circ s^n
 + \chi_{1,m} \cdot ( \chi_2 \circ s^n - \chi_{2,m} \circ s^n ) + \chi_{1,m} \cdot \chi_{2,m} \circ s^n  .       \]
Using the $s$-invariance of $\nu_\varphi$, we
get that $\left \vert  \int_{\Sigma} \chi_1 \cdot \chi_2 \circ s^n \,
 d\nu_\varphi \right \vert$ is less than or equal to 
\begin{equation*}\label{indp}
  \abs{\chi_1 -
 \chi_{1,m}}_1 \, \abs{\chi_{2}}_\infty  + \abs{\chi_1}_\infty \, \abs{\chi_2 -
 \chi_{2,m}}_1 + \big \vert {\int_{\Sigma} \chi_{1,m} \cdot \chi_{2,m} \circ s^n \, d\nu_\varphi } \big \vert .
\end{equation*}
Exponential estimates for $\abs{\chi_j -
 \chi_{j,m}}_1$ come from the $L^1(\nu_\varphi)$-cylinder assumption. We now focus on the last integral. Denoting $\chi_{j,m} := \sum_{C \in \AA_m} a_{j,C} 1_C$, we have 
\[  \big \vert {\int_{\Sigma} \chi_{1,m} \cdot \chi_{2,m} \circ s^n \, d\nu_\varphi } \big \vert \leq \sum_{C,C' \in \AA_m} \abs {a_{1,C}} \abs {a_{2,C'}} \abs{  \nu_\varphi(C \cap s^{-n}C') - \nu_\varphi (C) \nu_{\varphi}(C')  } . \]
Now we use $\abs {a_{j,C}} \leq \abs{\chi_j}_\infty$ and sum (\ref{ibbs}) over $(E,F) \in \AA_m \times s^{-n} \AA_m = \AA_{0,m} \times \AA_{n,m+n}$. By specifying
$m = [ n/2 ]$, we obtain $\left
\vert {\int_{\Sigma} \chi_{1,m} \cdot \chi_{2,m} \circ s^n \,
  d\nu_\varphi } \right \vert \leq c \abs{\chi_1}_\infty
\abs{\chi_2}_\infty e^{- \delta n /2}$.  \fin

The following result is due to Philipp-Stout (see \cite{PS},
Section 7).  Given an observable $\chi$, we set $S_n(\chi) := \sum_{j = 0}^{n-1} \chi \circ s^j$. 

\begin{thm}\label{philstout} Let $\chi : \Sigma \to \R$ be $\nu_\varphi$-centered and
  $L^p(\nu_\varphi)$-cylinder for some $p > 2$. Then we have:  
\begin{enumerate}
\item  ${1 \over \sqrt n} \abs{S_n(\chi)}_2$ converges to some $\si \geq 0$.
\item If $\si
> 0$, then $\chi$ satisfies the Almost Sure Invariance
Principle (ASIP).
\end{enumerate}
\end{thm}
The later means that there exist a sequence of random variables $(\SS_n)_n$ and  a
Brownian motion $\WW$ with variance $\si$ defined on a probability space
$(\Omega , \Prr)$ such that:

- $(\SS_0,\ldots, \SS_{n-1})$ and $(S_0(\chi) ,\ldots, S_{n-1}(\chi) )$ have the same distribution for $n \geq 1$,

-  there exists $c > 0$ such that $\SS_n =  \WW(n) + o(n^{1/2 -
  c})$ $\Prr$-almost surely. \\
  
 The proof of theorem \ref{philstout} consists in approximating $(S_n(\chi))_n$ by a martingale difference sequence. Property (\ref{ibbs}) provides a sufficient amount of mixing which ensures that approximation. We note that the ASIP implies several stochastic properties related to Brownian motion, like the Central Limit Theorem 	and the Law of Iterated Logarithm (see \cite{PS},
section 1).

%%%%%%%%%%%%%%%%%%%%%%%%%%%%%%%%%%%%%%%%%%%%%%%%%%%%%%%%%
\section{Proof of theorem A}\label{cocod} 
%%%%%%%%%%%%%%%%%%%%%%%%%%%%%%%%%%%%%%%%%%%%%%%%%%%%%%%%%

We obtain theorem A by using a coding tree technique. Let us recall the construction. Let $\CC_f$ denote the
critical set of $f$ and $\VV:= \cup_{i = 0}^\infty
f^i(\CC_f)$. Let $z \notin \VV$ and $\{ w_\al \, , \, \al \in \AA
\}$ be an enumeration of $f^{-1}(z)$. Assume that $(\ga_\al)_{\al \in \AA}$ is a collection of paths $[0,1] \to \Pj^k$ satisfying $\ga_\al(0) = z$, $\ga_\al(1) = w_\al$ and $\ga_\al [0,1] \cap \VV = \emptyset$ (theorem \ref{good} below will provide such paths). For $\tilde \al \in
\Sigma$ and $n \geq 1$ we define paths $\ga_n(\tal)$ and points
$z_n(\tal)$ as follows. First we set $\ga_1 (\tal) := \ga_{\al_0}$ and
$z_1 (\tal) := w_{\al_0}$. Now, assuming that $\ga_j (\tal)$ and $z_j
(\tal)$ have been defined for $1 \leq j \leq n$, we set
$\ga_{n+1}(\tilde \al)$ to be the lift of $\ga_{\al_n}$ by $f^{n+1}$ with starting point $z_n (\tal)$, and $z_{n+1}
(\tal) :=  \ga_{n+1} (\tal)(1)$. We observe that $\ga_n(\tal)$ and
$z_n(\tal)$ only depend on $[\al_0,\ldots,\al_{n-1}]$. We also note that $z_n : \AA_n \to f^{-n}(z)$ is a
bijection and that  $f \circ z_{n+1} = z_n \circ s$ holds for every $n \geq 1$.\\

The next result is theorem 3.2 of \cite{D1} with the bound $\theta_k$ specified (see below). That result allows to use the coding tree method and to prove the existence of coding maps $\om = \lim_n z_n$. Let $L_{z,w}$ denote the complex line in $\Pj^k$ joining $z$ and $w$.

\begin{thm}\label{good}
Let $\theta < \theta_k = {2 \over 5(k-1)}$. There exist a zero volume subset $\DD_\theta \subset \Pj^k$ and $\rho = \rho_\theta > 0$ satisfying the following properties. For every distinct points $(z,w) \in \Pj^k \setminus \DD_\theta \cup \VV$, there exist an injective smooth path $\ga :
[0,1] \to L_{z,w} \setminus \VV$ joining $(z,w)$, and decreasing topological discs $(\Delta_n)_n \subset L_{z,w}$ such that for every $n \geq n_{z,w}$:
\begin{enumerate}
\item $\ga[0,1] \subset \Delta_n$,
\item there exist $(1-d^{-\theta n }) d^{kn}$ inverse branches of $f^n$ on $\Delta_n$,
\item these branches satisfy $\Diam g_n (\Delta_n)  \leq  c \, d^{- \rho n }$.
\end{enumerate}
\end{thm}

That theorem precisely holds for $\theta$'s satisfying $\theta < 2\zeta$, $\theta + \zeta'' < 1$ and $(k-1)\theta < \zeta'' -\zeta'$ for some $0 < \zeta < \zeta' < \zeta'' <1$, those estimates come from section 3.2 of \cite{D1} (replace $d$ by $d^{k-1}$ in the definition of $\tau_{\theta_n}$ to fix a slip there). One can check that the value $\theta_k = {2 \over 5(k-1)}$ is consistent with the above conditions. \\

Let us now prove theorem A. For $\theta < \theta_k$ we set $\SS_\theta := \DD_\theta \cup f(\DD_\theta) \cup \VV
\cup \Per$, where $\DD_\theta$ is provided by theorem \ref{good} and $\Per$ is the set of periodic points. That set has zero volume. We now use the coding tree method. Given $z \notin \SS_\theta$ we fix an enumeration $f^{-1}(z) = \{ w_\al \}$ as before. We have $w_\al
\notin \VV \cup \DD_\theta \cup \{ z \}$ by definition of $\SS_\theta$. Let $(\ga_\al)_{\al \in \AA}$ be a collection of smooth paths
joining $(z, w_\al)$ given by theorem \ref{good}, and let $n_z:= \max n_{z,w_\al}$. Following the construction described above, theorem
\ref{good} asserts that for every $\al \in \AA$ and $n \geq
n_z$, there exist at least $(1-d^{-\theta n})d^{kn}$ cylinders
$[\al_0,\ldots,\al_{n-1}]$ satisfying:
\begin{equation}\label{decr}
\Diam  \ga_{n+1} [\al_0,\ldots,\al_{n-1},\al]  \leq  c \, d^{- \rho n }.
\end{equation}
We define $\BB_n \subset \AA_{n+1}$ as $\BB_n := \{ \Diam \ga_{n+1}[\al_0,\ldots,\al_n] > c \, d^{-\rho
  n} \}$. Let $\BB(n) := \bigcup_{p \geq n} \BB_p$ and
  $\GG(n) := \Sigma \setminus \BB(n)$. We recall that $d(.,.)$ denotes the standard distance on $\Pj^k$.
  
\begin{lem}\label{cdi}  For every $n \geq n_z$, we have:
\begin{enumerate}
\item $\cdb \BB_n \leq d^{k(n+1)} d^{-\theta n}$.

\item $d(z_n(\tal) ,   z_{n+1} (\tal)
  ) \leq c \,  d^{- \rho n}$ for every $\tal \in \GG(n)$.
\end{enumerate}
\end{lem}

\proof Let $\BB_n^\al$ be the set of
$(n+1)$-cylinders in $\BB_n$ whose last coordinate is equal to $\al$. Then (\ref{decr}) and the definition of $\BB_n$ yield $\cdb \BB_n^\al
\leq d^{-\theta n} d^{kn}$. Hence $\cdb \BB_n \leq d^{-\theta n}
d^{k(n+1)}$, which is the point 1. Now let $\tal \in \GG(n)$. The point 2 follows from $\tal \notin \BB_n$ and the observation $d(z_n(\tal),z_{n+1} (\tal) ) \leq \Diam \ga_{n+1}(\tal) = \Diam \ga_{n+1}[\al_0,\ldots,\al_n]$.  \fin

Using that $(\GG(n))_n$ is an increasing sequence of subsets, the point 2 of lemma \ref{cdi} yields $d(z_m(\tal) ,   z_{m+1} (\tal)
  ) \leq c \,  d^{- \rho m}$ for every $\tal \in \GG(n)$ and $m \geq n$. We deduce that $\om(\tal)
:= \lim_{m} z_m(\tal)$ exists for every  $\tal \in \GG(n)$ and that
\begin{equation}\label{starun}
\forall n \geq n_z \ , \ \forall \tal \in \GG(n) \ , \ d( z_n(\tal),
\om(\tal) ) \leq  c_\rho \, d^{-\rho n}.
\end{equation}
Now let $\varphi : \Sigma \to \R$ be an H\"older function such that $\tau_\theta (\varphi) = P_\varphi - \abs
{\varphi}_\infty - \log d^{k-\theta} > 0$ (we simply set $\tau = \tau_\theta (\varphi)$ in the sequel). 
\begin{lem}\label{jxx}
For every $n \geq n_z$, we have $\nu_\varphi(\BB_n) \leq  c \, e^{-n \tau}$. 
\end{lem}

\proof Using (\ref{ggiibb}) in section \ref{gibbs} and the definition of $\tau$, we obtain the following bound for the $\nu_\varphi$-measure of $(n+1)$-cylinders:
\begin{equation*}\label{cdeux}
\nu_\varphi [\al_0,\ldots,\al_n]   \leq  c_2  \,  e^{(n+1) (\abs{\varphi}_\infty - P_\varphi)} = c_2  \, d^{-(n+1)(k-\theta)} \, e^{-(n+1) \tau } .
\end{equation*}
We deduce $\nu_\varphi(\BB_n) \leq \cdb \BB_n \cdot c_2 \, 
d^{-(n+1)(k-\theta)} \, e^{-(n+1) \tau}$, which is less than $(c_2 \, d^\theta \, e^{-\tau}) e^{-n \tau}$ by lemma \ref{cdi}(1). \fin

We deduce from lemma \ref{jxx} and $\GG(n) = \Sigma \setminus \bigcup_{p \geq n} \BB_p$ that
\begin{equation}\label{stardeux}
\forall n \geq n_z \ , \ \nu_\varphi(\GG(n)) \geq 1 -  c_\tau  \, e^{- n
  \tau}.
\end{equation}
Now (\ref{starun}) and (\ref{stardeux}) immediately yield the points 1 and 2 of theorem A. Proposition \ref{exp} below yields the ASIP and the exponential decay of correlation (point 3 of theorem A) for H\"older observables $\psi : \Pj^k \to \R$. To see that, combine proposition \ref{exp} with theorems \ref{exdecay} and \ref{philstout} taking into account $\chi := \psi \circ \om$ and $f \circ \om
= \om \circ s$. 

\begin{prop}\label{exp}  Let $\psi : \Pj^k \to \R$ be a
  $\om_* \nu_\varphi$-centered H\"older observable. Then  $\chi := \psi \circ \om : \Sigma \to \R$ is a bounded $L^p(\nu_\varphi)$-cylinder observable for every $p \geq 1$.  
\end{prop}

\proof We obviously have $\chi \in L^p(\nu_\varphi)$, since $\psi$ is bounded. Now let us show that $\abs { \chi  - \E ( \chi \vert
  \AA_n  )  }_p \leq c \, e^{-n\ga}$ holds for some $\ga >0$. The proof consists in making an analysis on $\GG(n)$ and
$\BB(n) = \Sigma \setminus \GG(n)$  (see \cite{D1}, subsection 5.1.1).  For any Borel set
$A \subset \Sigma$, we set $\chi_A :=  \chi \cdot 1_A$. First observe that Jensen's inequality
and (\ref{stardeux}) yield  
\[ \abs { \chi_{\BB(n)}   - \E ( \chi_{\BB(n)} \vert \AA_n  )
  }_p \leq 2 \abs{\chi}_\infty  \, \nu_\varphi (\BB(n))^{1/p}   \leq 2 \abs{\chi}_\infty  \, \left( c_\tau   \, e^{- n  \tau}
\right)^{1/p}. \] 
It remains to prove the following estimate for some $\lambda > 0$:
\begin{equation}\label{esb2}
\abs { \chi_{\GG(n)} - \E ( \chi_{\GG(n)} \vert \AA_n  )
  }_p \leq c \, e^{-n \lambda}.
\end{equation}
We set $\phi:= \chi_{\GG(n)} - \E ( \chi_{\GG(n)} \vert
\AA_n)$ and estimate $\abs
   {\phi_{\BB(n)}}_p$ and
   $\abs{\phi_{\GG(n)}}_p$. Observing that $\phi_{\BB(n)} =  - \E ( \chi_{\GG(n)} \vert
\AA_n)$ on $\BB(n)$ and zero elsewhere, we deduce from (\ref{stardeux}):
\begin{equation}\label{raz}
 \abs{\phi_{\BB(n)}}_p \leq \abs { \E ( \chi_{\GG(n)} \vert
\AA_n)  }_{2p} \cdot \nu(\BB(n)) ^{1/2p} \leq \abs{\chi}_{2p} \cdot \left( c_\tau \, e^{- n  \tau} \right)^{1/2p}.
\end{equation}
Concerning $\abs{\phi_{\GG(n)}}_p$, we have for every $\tal \in
\GG(n)$:
\begin{equation*}\label{lki}
  \phi_{\GG(n)}(\tal) = \int_{[\tal]_n \cap \GG(n)} \left( \chi(\tal)
-\chi (\tbe) \right)  d\nu_\tal(\tbe)  + \chi(\tal) \cdot  \nu_\tal([\tal]_n \cap \BB(n)),
\end{equation*}
where $[\tal]_n := [\al_0,\ldots,\al_{n-1}]$ and $\nu_\tal$ is the conditional
measure of $\nu_\varphi$ on $[\tal]_n$. We deduce from $\chi =  \psi \circ \om $, (\ref{starun}) and the
fact that  $\psi$ is H\"older (say of exponent $h$):
\begin{equation*}
 \forall \tal \in \GG(n) \ , \ \abs{\phi_{\GG(n)}(\tal)} \leq \left( 2 \, 
 c_\rho \, d^{-\rho n} \right)^h  + \abs{\chi}_\infty \cdot  \nu_\tal([\tal]_n \cap \BB(n)).
\end{equation*}
Integrating over $\GG(n)$ and using (\ref{stardeux}), we obtain $\abs{\phi_{\GG(n)}}_p^p \leq d^{- \rho n h p} +  \abs{\chi}^p_\infty  \cdot  c_\tau \, e^{- n
 \tau}$ up to a multiplicative constant. Combining (\ref{raz}) with that estimate, we get (\ref{esb2}) as desired. \finsec

%%%%%%%%%%%%%%%%%%%%%%%%%%%%%%%%%%%%%%%%%%%%%%%%%%%%%%%%%
\section{Proof of theorem B} \label{P1}
%%%%%%%%%%%%%%%%%%%%%%%%%%%%%%%%%%%%%%%%%%%%%%%%%%%%%%%%%

Given $\theta < \theta_k$, theorem A yields a zero volume subset $\SS_\theta$  providing coding maps $\om : \Sigma \to \Pj^k$ and invariant measures $\nu = \om_* \nu_\varphi$ on $\Pj^k$. Now our aim is to show the

\begin{theorem}{ B:}
There exists a zero volume subset $\HH \subset \Pj^k$ satisfying the
  following property. Let $\theta < \theta_k$ and $z \notin \SS_\theta \cup \HH$. Let $\varphi : \Sigma \to \R$ be an H\"older function satisfying
  $\tau_\theta(\varphi) > 0$ and $\om : \Sigma \to \Pj^k$  provided by theorem A. Then $h_{\om_* \nu_\varphi} (f) = h_{\nu_\varphi}(s)$.
\end{theorem}

 We shall use the following Abramov-Rohlin's formula \cite{AR}. We refer also to the article of Ledrappier-Walters \cite{LW}, lemma 3.1. The statement written here is adapted to our context, the results in \cite{AR} and \cite{LW} are more general.
  
 \begin{thm}\label{aroh}
 Let $\om : \Sigma \to \Pj^k$ be a measurable map satisfying $f \circ
 \om = \om \circ s$ and $\nu_\varphi$ be a Gibbs measure on
 $\Sigma$. Then the $f$-invariant measure $\nu = \om_*\nu_\varphi$ satisfies:  
  \[ h_{\nu_\varphi}(s)= h_{\nu}(f) +  \lim_{n \to \infty} {1 \over n}
  \, H_{\nu_\varphi} (\AA_n \wert \om^{-1}\epsi) . \] 
 \end{thm}

Let us specify the notations. Given a measurable partition $\xi$, we denote $H_{\nu_\varphi} (\AA_n
\wert \xi) = \int_{\Sigma} \, H_{\xi(\tal)} (\AA_n) \, d \nu_\varphi
(\tal)$, this is the conditional entropy of $\AA_n$ with respect to $\xi$. We recall that $H_{\xi(\tal)} (\AA_n) = - \sum_{
C \in \AA_n}   \log
\nu_{\xi(\tal)} (C) \, \nu_{\xi(\tal)}(C)$, where $\xi(\tal)$ is the atom of $\xi$ containing
$\tal$ and $\nu_{\xi(\tal)}$ is the conditional
measure of $\nu_\varphi$ on  $\xi(\tal)$. Finally, in the statement of theorem \ref{aroh}, $\epsi$ stands for the partition of $\Pj^k$ into points. \\

In view of theorem \ref{aroh}, theorem B can  be reformulated as follows: 

\begin{theorem}{ B':} There exists a zero volume subset $\HH \subset \Pj^k$ satisfying the following property. Let  $z \notin  \SS_\theta \cup  \HH$ and $\varphi : \Sigma \to \R$ be an H\"older function satisfying
  $\tau_\theta(\varphi) > 0$. Let $\om : \Sigma \to \Pj^k$  provided by theorem A. Then $\lim_{n \to \infty} {1 \over n} \, H_{\nu_\varphi}(\AA_n \wert \om^{-1}\epsi) = 0$. 
\end{theorem}

We follow for the proof the approach of Przytycki  \cite{P} employed for $\om$ boundary extensions of Riemann mappings defined on the unit disc in $\C$.

%%%%%%%%%%%%%%%%%%%%%%%%%%%%%%%%%%%%%%%%%%%%%%%%%%%%%%%%%
\subsection{Preliminaries}\label{prel} 
%%%%%%%%%%%%%%%%%%%%%%%%%%%%%%%%%%%%%%%%%%%%%%%%%%%%%%%%%

We deal in this subsection with lemmas related to Misiurewicz-Przytycki's inequality \cite{MP}. We shall work with an iterate $g : = f^q$. The first lemma concerns the tree of preimages of points. Let $w \in \Pj^k$ and
$G_m \subset g^{-m}(w)$. Our aim is to estimate the cardinal of $G_m$ in terms of a branching condition. For $1
\leq j \leq m$ we set $G_{m-j}  := g^j (G_m)$ and say that $p \in G_{m-j}$ is \emph{branching} if $\cdb
g^{-1}(p) \cap G_{m-j+1}   \geq 2$. Let also
\[ \forall y \in G_m \ , \  \TT_m(y) := \{  \, 0 \leq j \leq m-1 \, , \, g^{j+1} (y) \textrm{ is branching}  \, \} .\] 
  
\begin{lem}\label{trr}
If $\cdb \TT_m(y) \leq s$ for every $y \in G_m$, then $\cdb G_m \leq d^{qks}$.  
\end{lem}

\proof  We proceed by induction on $m$. The lemma is clear when $m=1$,
since $s \leq 1$ in that case. Assume now that the assertion holds for some $m \geq 1$. Let $G_{m+1}
\subset g^{-(m+1)}(w)$ satisfying $\cdb \TT_{m+1}(y) \leq s$ for every
$y \in G_{m+1}$. We distinguish two cases. If $w$ is not branching, then $\cdb
\TT_{m}(y) =  \cdb \TT_{m+1}(y) \leq s$ for every $y \in G_{m+1}$, hence $\cdb G_{m+1} \leq d^{qks}$ by induction. If $w$ is branching, let
$\{ w_1, \ldots , w_r \} := g^{-1}(w) \cap G_1$ (with $2 \leq r \leq
d^{qk}$) and $G_m^i := \{ y \in G_{m+1} \, , \,  g^m(y) = w_i
\}$. Observing that $\cdb \TT_{m}(y) \leq s-1$ for every $y \in
G_m^i$ and using the induction, we get $\cdb G_m^i \leq d^{qk(s-1)}$. That implies $\cdb G_{m+1} \leq r d^{qk(s-1)} \leq d^{qks}$ as desired. \fin

We shall need two other lemmas. Let $\Jac g$ be the smooth function $\Pj^k \to \R^+$ satisfying
$g^*\om^k = \Jac g \cdot \om^k$, where $\om^k$ denotes the standard volume
form of $\Pj^k$. Note that $\{ \Jac g = 0 \}$ coincides with the
critical set $\CC_g$ of $g$. For every $\delta > 0$, let $\CC_g(\delta) := \{ \Jac g \leq \delta \}$. We set $L := \abs{\Jac g}_\infty + 1$, $L' := \abs { d \, \Jac
 g}_\infty + 1$, $a < 1$ and $\ga < 1$.\\

Let $\delta_\ga := \left( {a \over L} \right)
 ^{1/\ga}$, and for every $y \in \Pj^k$,
 \[  \HH_m(y) := \left \{ \, 0 \leq i \leq m-1 \, , \, g^i(y) \in
 \CC_g(\delta_\ga) \, \right\}  .\]
 Let $\HH_m[\ga] := \left \{  \, y \in \Pj^k \, , \,   \cd \, \HH_m(y)
 > m \ga  \,  \right \}$ and $\HH[\ga] := \limsup_m g^m(\HH_m[\ga])$. We shall write $\HH^q[\ga]$ for $\HH[\ga]$ to insist on the dependence on $q$.
 
\begin{lem}\label{fiber1} We have $\Vol \ \HH^q[\ga]  = 0$.
\end{lem}

\proof For every $y \in
\HH_m[\ga]$, we get $\Jac g^m (y) \leq (\delta_\ga) ^{\ga m}
L^{(1-\ga)m}  \leq (\delta_\ga) ^{\ga m} L^m = a^m$. That implies
$\Vol \ g^m(\HH_m[\ga]) \leq  a^m \, \Vol \
\HH_m[\ga]  \leq a^m$. Hence $\Vol \ \HH^q[\ga] \leq c \,  a^m$ for every $m \geq 1$, that completes the proof of the lemma. \fin 

The next lemma follows from the compacity of $\Pj^k$:

\begin{lem}\label{fifiber}  There exists $\kappa(\ga) = \kappa(\ga,q) < {\delta_\ga} / 2L'$ satisfying: 
\begin{equation*}\label{sep}
 \forall (y , y') \in \Pj^k \setminus \CC_g(\delta_\ga /2) \ , \   d(y,y') \leq \kappa(\ga) \,  \textrm{ and } \, y \neq y' \, \Longrightarrow \, g(y) \neq g(y'). 
\end{equation*}
\end{lem}

Using that definition of $\kappa(\ga)$, we define for every $y \in \Pj^k$:   
\[ \FF_m(y) :=  \left \{   0 \leq i \leq m-1 \, , \, \exists h \in \Pj^k \, , \,
d(g^i(y),h) \leq \kappa(\ga) \, , \,   g^i(y) \neq h    \textrm{ and }
g^{i+1}(y) = g(h) \right \} . \]

\begin{lem}\label{fiber} For every $z \notin \HH^q[\ga]$ there exists $m(\ga) = m(\ga,q,z)\geq 1$ such that 
\[ \forall m \geq m(\ga) \ , \ \forall y \in g^{-m}(z) \ , \ \cdb \FF_m(y) \leq \ga m . \]  
\end{lem}

\proof  Let $m(\ga) \geq 1$ be such that $z \notin g^m(\HH_m[\ga])$
for every $m \geq m(\ga)$. By definition of $\HH_m[\ga]$, we have $\cdb \HH_m(y) \leq \ga m$  for every $m \geq m(\ga)$ and 
$y \in g^{-m}(z)$. Hence it suffices to prove $\FF_m(y) \subset
\HH_m(y)$ to complete the proof of the lemma. 

Let $i \in \FF_m(y)$. From the definition of $\kappa(\ga)$,
either $g^i(y)$ or $h$ is in
$\CC_g(\delta_\ga/2)$. That implies that $g^i(y)$ and $h$ are both in
$\CC_g(\delta_\ga)$, because $\abs{ \Jac g( g^i(y) ) - \Jac g ( h ) } \leq L' d
(g^i(y),h) \leq L' \kappa(\ga) \leq {\delta_\ga} /2$. Hence $i \in
\HH_m(y)$ as desired.  \finsec

%%%%%%%%%%%%%%%%%%%%%%%%%%%%%%%%%%%%%%%%%%%%%%%%%%%%%%%%%
\subsection{Proof of theorem B'} 
%%%%%%%%%%%%%%%%%%%%%%%%%%%%%%%%%%%%%%%%%%%%%%%%%%%%%%%%%

Let $(\ga_p)_{p \geq 1}$ satisfying $\lim_p \ga_p = 0$ and let $\HH := \cup_{p,q \geq 1} \HH^q[\ga_p]$. That subset has zero volume
by lemma \ref{fiber1}. Let $\theta < \theta_k$ and fix $z \notin   \SS_\theta \cup  \HH$. Let $\varphi : \Sigma \to \R$ be an H\"older function satisfying
  $\tau_\theta(\varphi) > 0$ and $\om = \lim_n z_n : \Sigma \to \Pj^k$ be a coding map given by theorem A. Let also $n_z$ provided by that theorem. We establish in this subsection: 
\begin{equation}\label{larss}
\forall p, q \geq 1 \ , \ \limsup_{n \to + \infty} {1 \over n} \, H_{\nu_\varphi}(\AA_n \wert \om^{-1}\epsi)
\leq 6 \ga_p \log d^k + {1 \over q} \log 2. 
\end{equation}
That implies $\lim_n {1 \over n} \, H_{\nu_\varphi}(\AA_n
\wert \om^{-1}\epsi) = 0$ as desired. \\

Let us fix $p,q \geq 1$. We consider $\kappa(\ga_p)$ and $m(\ga_p)$ given by lemmas \ref{fifiber} and \ref{fiber} (they depend on $q$ and $z$). Let $n_1 \geq n_z$ be such that (see (\ref{starun})
and (\ref{stardeux})):
\begin{equation}\label{vco} 
 \nu_\varphi(\GG(n_1)) > 1 - \ga_p \ \textrm{ and } \  \forall n \geq n_1 \ , \ \forall \tal \in \GG(n_1)   \ , \  d(z_n(\tal) , \om(\tal)) < \kappa(\ga_p) / 4 . 
 \end{equation}
By Birkhoff's theorem, there exist $E \subset \Sigma$, $m_1 \geq 1$ such that $\nu_\varphi(E) > 1 - \ga_p$ and 
\begin{equation}\label{pmmk} 
\forall \tal \in E \ , \ \forall m \geq m_1 \, , \,  \cdb \{  \, 0 \leq i \leq m-1 \, , \, s^{qi} (\tal) \in
\GG(n_1)   \, \}  > (1 - 2\ga_p)m . 
\end{equation}
We set $n_{p,q} := n_1 + q m_1 + q  m(\ga_p)$. \\

 We introduce the partitions $\QQ := \{  E , E^c  \}$
and $\PP := \om^{-1}\epsi \vee  \QQ$.  The next lemma allows to replace
$H(\AA_n \wert  \om^{-1}\epsi)$ by $H(\AA_n \wert  \PP)$ for the proof
of (\ref{larss}) (we  denote
$H$ for $H_{\nu_\varphi}$).

\begin{lem}\label{low}
$H(\AA_n \wert  \om^{-1}\epsi) \leq H ( \QQ  ) + H (\AA_n  \wert \PP )$.
\end{lem}

\proof We use twice $H(\zeta \vee \xi \wert \xi') = H(\zeta \wert \xi')
+ H( \xi  \wert \xi' \vee \zeta )$ (\cite{KH}, section 4.3) to get:
\[ H ( \AA_n \vee  \QQ \wert \om^{-1}\epsi )   =     H (  \AA_n \wert \om^{-1}\epsi  )  + H  (\QQ  \wert \om^{-1}\epsi \vee  \AA_n   )  
                                            =     H ( \QQ  \wert \om^{-1}\epsi ) + H (\AA_n \wert \PP ). \]
We deduce $H (  \AA_n \wert \om^{-1}\epsi  ) \leq  H ( \QQ \wert
\om^{-1}\epsi ) +  H (\AA_n \wert \PP ) \leq H ( \QQ  ) + H (\AA_n  \wert \PP )$.\fin 

We now write 
\begin{equation}\label{hen}
 H (\AA_n \wert \PP )  =   \int_{E}   \, H_{\PP(\tal)} (\AA_n) \, d
 \nu_\varphi (\tal) + \int_{E^c}   \, H_{\PP(\tal)} (\AA_n) \, d
 \nu_\varphi (\tal) 
\end{equation}
and estimate those integrals. 

\begin{lem}\label{hen1}
For every $n \geq 0$, $\int_{E^c} \, H_{\PP(\tal)} (\AA_n) \, d
\nu_\varphi (\tal) \leq \ga_p  \log d^{kn}$.
\end{lem}
The proof comes from $H_{\PP(\tal)} (\AA_n) \leq \log \cdb \AA_n = \log d^{kn}$ and $\nu_\varphi(E^c)\leq \ga_p$. The integral over $E$ is more delicate. We shall prove in next subsection:

\begin{lem} \label{hen2}
For every $n \geq n_{p,q}$, $\int_{E} \, H_{\PP(\tal)} (\AA_n) \, d
\nu_\varphi (\tal) \leq \log d^{k(5{\ga_p} n + 2n_1 + q)} + {n \over q} \log 2$.
\end{lem}
For now, the three previous lemmas and (\ref{hen}) imply:
\[  \forall n \geq n_{p,q} \ , \  H (\AA_n \wert \om^{-1}\epsi ) \leq
H ( \QQ  ) +   \ga_p  \log d^{kn} + \log d^{k(5{\ga_p} n + 2 n_1 + q)}  + {n \over q} \log 2  . \]
That yields $\limsup_{n \to \infty} {1 \over n} \, H(\AA_n \wert
\om^{-1}\epsi) \leq  6 {\ga_p} \log d^k + {1 \over q} \log 2$, which
is (\ref{larss}).

%%%%%%%%%%%%%%%%%%%%%%%%%%%%%%%%%%%%%%%%%%%%%%%%%%%%%%%%%
\subsection{Proof of  lemma \ref{hen2}} \label{P2}
%%%%%%%%%%%%%%%%%%%%%%%%%%%%%%%%%%%%%%%%%%%%%%%%%%%%%%%%%

We recall that $\om = \lim_n z_n$ and $\PP =
\om^{-1}\epsi \vee  \{  E , E^c \}$. The proof of lemma \ref{hen2} consists in studying the size of the fibers of $\om$. For every $x \in
\Pj^k$, we set $E(x) :=  \om^{-1} \{ x \}  \cap E$. Let $\AA_n(x)  := \pi_n  \big (
E(x)  \big)$, where $\pi_n : \Sigma \to \AA_n$ denotes the projection $\pi_n(\tal) := [\al_0, \ldots , \al_{n-1}]$. For every $C \in
\AA_{n_1}$ and $n \geq n_1$, we set $\AA_n^C(x) := \AA_n(x) \cap \pi_{n_1}^{-1} (C)$.\\

We claim that lemma \ref{hen2} is a consequence of the following proposition.

\begin{prop}\label{SN}
$\forall x \in \Pj^k \, , \, \forall n \geq n_{p,q}  \, , \, \cdb
  \AA_n(x) \leq  d^{k(5{\ga_p} n + 2 n_1 + q)} \cdot 2^{n/q}$.
  \end{prop}
 Observe indeed that $\cdb \AA_n(x)$ is the cardinal of the
 partition of $E(x)$ induced by the set of $n$-cylinders $\AA_n$. Proposition \ref{SN} therefore implies $H_{E(x)} (\AA_n) \leq \log d^{k(5 {\ga_p} n +  2 n_1 + q )} + {n
 \over q} \log 2$. Lemma \ref{hen2} then follows from (note that $\PP(\tal) = E(x)$ for every $\tal \in E(x)$):
\begin{equation*} 
 \int_{E} \, H_{\PP(\tal)} (\AA_n) \, d \nu_\varphi (\tal) = \int_{\Pj^k} \, \left[ \, \int_{E(x)} \, H_{E(x)} (\AA_n) \, d \nu_{E(x)} (\tal) \, \right] \, d (\om_*\nu_\varphi) (x)
\end{equation*}
and the fact that $\nu_{E(x)}$ and $\om_*\nu_\varphi$ are probability
measures  on $E(x)$ and $\Pj^k$ respectively. \\

Now we show proposition \ref{SN}. We shall need  lemmas \ref{hhh} and \ref{hgt} below. We recall that $g = f^q$ and
  $n_{p,q} = n_1 + q m_1 + q  m(\ga_p)$. For every $n \geq
  n_{p,q}$ we set $m := [(n-n_1)/q]$. We therefore have $n = m q + n_1 + r$ for some $0 \leq r \leq
  q-1$, and $m \geq \max \{ m_1 , m(\ga_p) \}$.  We fix $x \in \Pj^k$ and define for every $(a,b) \in \AA_n(x) = \pi_n(E(x))$: 
\[ \LL_m(a,b) := \{ \,  0 \leq i \leq m-1 \, , \, d( g^i ( z_n (a) )  , g^i ( z_n (b) )  ) > \kappa(\ga_p) / 2  \,  \}.\]
The next lemma is established in subsection \ref{prooo}. The introduction of the iterates of $f$ (namely the definition of $E$, see (\ref{pmmk})) will be crucial for the proof.
\begin{lem}\label{hhh}
For every $n \geq n_{p,q}$ and $(a,b)  \in \AA_n(x)$, $\cdb \LL_m(a,b) \leq 4 \ga_p m$. 
\end{lem}

We fix for the sequel $C \in \AA_{n_1}$ and $a = [a_0, \ldots , a_{n-1}] \in \AA_n^C(x)$ (if not
empty). We also fix $\Delta \subset \{ 0 , \ldots , m-1 \}$ satisfying $\cdb
\Delta \leq 4 \ga_p m$. Note that by setting 
\[ \AA_n^C(x,\Delta) := \big \{ \ b  \in \AA_n^C(x) \, , \, \LL_m(a,b)  \subset \Delta \ \big \} ,  \]
we have $\AA_n^C(x) = \bigcup_{\cdb \Delta \leq 4 \ga_p m} \, \AA_n^C(x,\Delta)$ (see lemma \ref{hhh}). That decomposition will be useful at the end of this subsection. \\

Let us observe that $z_n(\AA_n^C(x,\Delta)) \subset f^{-n}(z)$ and $n=mq + n_1 + r$. That implies  
\[ g^m (   z_n  (  \AA_n^C(x,\Delta) )   ) \subset  f^{-(n_1+r)}(z). \]
Given $w \in f^{-(n_1+r)}(z)$, we define $G_m(w) := z_n(\AA_n^C(x,\Delta))  \cap g^{-m}(w)$. Our aim is to prove $\cdb G_m(w)  \leq d^{k \cdot 5
  \ga_p n}$. \\
  
We apply for that purpose lemma \ref{trr} with $G_m = G_m(w)$. Let us recall the definitions: $G_{m-i}(w)  := g^i ( G_m (w))$ for $1 \leq i \leq m$ and we say that $p \in G_{m-i}(w)$
is branching if $\cdb g^{-1}(p) \cap G_{m-i+1}(w)  \geq 2$. For $y \in G_m(w)$, we define the subsets of $\{  0 , \ldots , m-1\}$:
 \[ \FF_m(y) =  \left \{    \exists h \in \Pj^k \, , \,
d(g^i(y),h) \leq  \kappa(\ga_p) \, , \,   g^i(y) \neq h    \textrm{ and }
g^{i+1}(y) = g(h) \right \}  , \]
\[ \TT_m(y) = \{   g^{i+1}(y) \textrm{ is
  branching} \} . \] 
The next lemma will be proved in subsection \ref{prooo}. 
\begin{lem}\label{hgt}
For $n \geq n_{p,q}$ and $y \in G_m(w)$, we have $\TT_m(y) \subset \Delta \cup \FF_m(y)$.
\end{lem} 
Now let us deduce the desired estimate for $\cdb G_m(w)$. Since $z \notin \HH^q(\ga_p)$ and $m \geq m(\ga_p)$, lemma \ref{fiber} yields $\cdb \FF_m(y) \leq {\ga_p} m$. We therefore obtain $\cdb \TT_m(y) \leq  5 \ga_p m$ from lemma  \ref{hgt}. Finally, by applying lemma \ref{trr} with $s=5\ga_p m$, we get $\cdb G_m(w) \leq d^{qk \cdot 5 \ga_p m} \leq d^{k \cdot 5   \ga_p n}$. \\

Let us complete the proof of proposition \ref{SN}.  From the observation 
  \[ \AA_n^C(x,\Delta) = \bigcup_{w \in f^{-(n_1+r)}(z)} G_m(w) \,  , \]
  we get $\cdb  \AA_n^C(x,\Delta)    \leq d^{k \cdot (5
  \ga_p n + n_1 + q)}$. Hence $\cdb \AA_n(x) \leq  d^{k \cdot (5 \ga_p n + n_1 + q )}  d^{k
  n_1} 2^m$ by taking different $C$ and $\Delta$. The desired upper bound follows using $m \leq n / q$.

%%%%%%%%%%%%%%%%%%%%%%%%%%%%%%%%%%%%%%%%%%%%%%%%%%%%%%%%%
\subsection{Proof of lemmas \ref{hhh} and \ref{hgt}}\label{prooo}
%%%%%%%%%%%%%%%%%%%%%%%%%%%%%%%%%%%%%%%%%%%%%%%%%%%%%%%%%

Let us begin with the \textsc{proof of lemma \ref{hhh}}. Let $n \geq n_{p,q}$ and $(a,b) \in \AA_n(x)=\pi_n(E(x))$. We want to prove $\cdb \LL_m(a,b) \leq 4 \ga_p m$. Assume to the contrary that $\cdb
\LL_m(a,b) > 4 \ga_p m$. Let $(\tal, \tbe) \in E(x)$ such that
$\pi_n(\tal) = a$ and $\pi_n(\tbe) = b$. Using (\ref{pmmk}), the definition of $\LL_m(a,b)$ and  $\cdb
\LL_m(a,b) > 4 \ga_p m$, there exists
$0 \leq i \leq m-1$ satisfying:
\begin{equation}\label{znh} 
 \delta_n := d( g^i ( z_n (\tal) )  , g^i ( z_n (\tbe) )  ) > \kappa(\ga_p) / 2  \ \textrm{ and }  \ s^{qi}(\tal) , s^{qi}(\tbe) \in \GG(n_1).    
 \end{equation}
Now let us fix $l \geq 1$. We have (use $g^i \circ z_n = z_{n-qi} \circ
s^{qi}$ for the second line):
\[   \delta_n \leq d( g^i ( z_n (\tal) )  , g^i ( z_{l+qi} (\tal) )  ) + d( g^i (
z_{l+qi} (\tal) )  , g^i ( z_{l+qi} (\tbe) )  ) + d( g^i ( z_{l+qi}
(\tbe) )  , g^i ( z_n (\tbe) )  )   \]
\[   = d( z_{n-qi} (  s^{qi}(\tal) )  ,  z_l (s^{qi } (\tal) )  ) + d(
g^i ( z_{l+qi} (\tal) )  , g^i ( z_{l+qi} (\tbe) )  ) + d( z_l
(s^{qi}(\tbe) )  ,  z_{n-qi} ( s^{qi} (\tbe) )  )  . \]
Now letting $l \to + \infty$, we obtain
\[  \delta_n \leq d( z_{n-qi} (  s^{qi}(\tal) )  ,  \om (s^{qi } (\tal) )  ) + d( g^i ( \om (\tal) )  , g^i ( \om (\tbe) )  ) + d( \om (s^{qi}(\tbe) )  ,  z_{n-qi} ( s^{qi} (\tbe) )  )  . \]
But the middle term vanishes since $\om(\tal) = \om(\tbe)$ equal to $x$. Hence
\[  \delta_n \leq  d( z_{n-qi} (  s^{qi}(\tal) )  ,  \om (s^{qi } (\tal) )  ) + d( \om (s^{qi}(\tbe) )  ,  z_{n-qi} ( s^{qi} (\tbe) )  )  .  \]
We now use (\ref{vco}) with $s^{qi}(\tal) , s^{qi}(\tbe) \in \GG(n_1)$ (note that $n-qi \geq n-qm \geq  n_1$) to get
\[ \delta_n \leq  \kappa(\ga_p) / 4 + \kappa(\ga_p) /4 = \kappa(\ga_p) / 2. \] 
That contradicts the first part of (\ref{znh}) and proves lemma \ref{hhh}. \\

We now deal with the \textsc{proof of lemma \ref{hgt}}. Let $w \in f^{-(n_1+r)}(z)$ and $y \in G_m(w)$. Let also 
$i \in \TT_m(y)$: that means that $g^{i+1}(y)$ is  branching. Hence there exists $y' \in G_m(w)$
such that $g^i(y) \neq  g^i(y')$  and $g^{i+1}(y) = g^{i+1}(y')$. Let $(b,b') \in \AA_n^C(x,\Delta)$ such that  $y = z_n(b)$ and $y' =
z_n(b')$. Now we assume that $i \notin \Delta$ and we prove that $i \in \FF_m(y)$. Let us write (recall that some $a \in \AA_n^C(x)$ has been fixed at this stage): 
\[ d( g^i (y) , g^i (y') ) \leq   d \left( g^i(z_n (b)) , g^i (z_n(a)) \right)  + d \left ( g^i(z_n(a)) , g^i (z_n (b')) \right) . \]
Using $\LL_m(a,b)  , \LL_m(a,b')  \subset \Delta$ (from $(b,b') \in \AA_n^C(x,\Delta)$) and $i \notin \Delta$, we obtain that $d( g^i (y) , g^i (y') ) \leq   \kappa(\ga_p)$. We finally get $i \in \FF_m(y)$ by taking $h = g^i(y')$ in the definition of $\FF_m(y)$. That completes the proof of lemma \ref{hgt}.

%%%%%%%%%%%%%%%%%%%
\section{Proofs of theorems C and D}\label{cetd}
%%%%%%%%%%%%%%%%%%%

We prove in this section the following result (we set $x^+ := \max \{  x, 0\}$). 

\begin{theorem}{ D:}  Let $f$ be an holomorphic endomorphism of $\Pj^k$ of degree $d \geq 2$ and $m$ be an ergodic
$f$-invariant measure. Let $\Lambda_1 > \ldots >
\Lambda_q \geq - \infty$ denote the distinct Lyapunov exponents of $m$ and $(m_i)_{1 \leq i \leq q}$  their multiplicities. We assume that $q \geq 2$. Then for every $2 \leq j \leq q$, the metric entropy of $m$ satisfies
\[  h_m  \leq \log d^{m_1 + \ldots + m_{j-1}}  + 2 m_j \Lambda_j^+ +
\ldots + 2 m_q \Lambda_q^+ .\]
\end{theorem}

That result extends Margulis-Ruelle's inequality $h_m \leq  2
m_1 \Lambda_1^+ + \ldots + 2 m_q \Lambda_q^+$. We recall that  theorem D was proved by de Th\'elin \cite{dT}
assuming $\Lambda_q > -\infty$. Here we extend it to $\Lambda_q \geq -
\infty$. Let us see before how theorem D implies the

\begin{theorem}{ C:}
If $h_m > \log d^{k-1}$, then $\Lambda_1 > \ldots >
\Lambda_q \geq {1 \over 2} (h_m - \log d^{k-1} )$.
\end{theorem}

When $m$ has $q \geq 2$ different exponents, theorem C follows from the inequality  $\Lambda_q^+ \geq
  {1 \over 2 m_q}(h_m - \log d^{k-m_q})$ given by theorem D, and from the observation
  \[ (1-{1 \over  m_q})h_m \leq  (1-{1 \over  m_q})\log d^k = \log
  d^{k-1} - {1 \over m_q}\log d^{k-m_q} , \]
  which implies ${1 \over  m_q} (h_m - \log d^{k-m_q}) \geq h_m - \log d^{k-1}$. When $m$ has a single exponent $\Lambda$, we use the inequality $\Lambda^+ \geq {1 \over 2k} h_m$, given by Margulis-Ruelle's
  inequality, and the observation $(1-{1 \over  k}) h_m \leq  \log
  d^{k-1}$ (take the preceding estimate with $m_q = k$).

%%%%%%%%%%%%%%%%%%%
\subsection{Preliminaries}\label{ospe}
%%%%%%%%%%%%%%%%%%%

Let us set $\Lambda_q = - \infty$ to fix the ideas. We can
  assume that $\Lambda_1 > 0$, otherwise $h_m = 0$ by Margulis-Ruelle's inequality, and theorem D is obvious. We fix $1 \leq u \leq q-1$ such that $\Lambda_1 >
  \ldots > \Lambda_u > 0 \geq \Lambda_{u+1} > \ldots >
  \Lambda_q$. It clearly suffices to establish theorem D for $2 \leq j \leq u+1$. Let $p := m_{u+1} + \ldots +
m_q$ and $q_j := m_{j} + \ldots + m_u$ for $1 \leq j \leq u$. \\
   
Let  $(\tau_x)_{x \in \Pj^k}$ be a family of charts satisfying
   $\tau_x : \D^k(\rho_0) \to \Pj^k$, $\tau_x(0)=x$ and such that $\tau_x^{\pm 1}$
   have bounded derivatives. Let $f_x :=
\tau_{f(x)}^{-1} \circ f \circ \tau_x$ and $f_x^n := f_{f^{n-1}(x)}
\circ \ldots \circ f_x$. These mappings are defined on a neighbourhood of the origin in $\C^k$. 

Let $\orb := \{ \hat x = (x_n)_{n \in \Z} \,
   , \, x_{n+1} = f(x_n) \}$ be the set of orbits and $\hat f$ be the left shift acting on $\orb$. We let $\pi(\hat x) :=
   x_0$ and $\hat m$ be the $\hat f$-invariant measure on $\orb$ satisfying $\hat m (\pi^{-1}A) = m(A)$ for every borel set $A \subset \Pj^k$.  We set $\hat x_n :=
   {\hat f}^n (\hat x)$ for every $n \in \Z$ and $\hat x_0 = \hat x$. We say that $\varphi_\epsi : \orb \to [1, +\infty[$ is $\epsi$-tempered if $e^{-\epsi} \varphi_\epsi
    \leq \varphi_\epsi \circ \hat f  \leq   e^{\epsi} \varphi_\epsi$.  \\

  The following result is due to Newhouse (see \cite{N1}, theorem
2.3). It is adapted to our context since it provides an Oseledec-Pesin's theorem in the case $\Lambda_q = -
\infty$. The statement we give here focus on the non
negative part of the spectrum: analog properties hold for the negative
part, but we shall not need them. For a linear subspace $E \subset
\C^k$, we denote $E^* := E \setminus \{ 0 \}$. In what
follows, $i$ ranges $\{ 1, \ldots , u \}$.
    
\begin{thm}\label{newhouse}  For $\hat m
  \pp$ $\hat x$, there exist splittings $\C^k
  = \oplus_{i=1}^{u}  E_i(\hat x) \oplus E_{u+1}(\hat x)$ such that:
\begin{enumerate}
\item $d_0 f_x : E_i(\hat x) \to E_i(\hat x_1)$ is invertible and $d_0 f_x (E_{u+1}(\hat x)) \subset E_{u+1}(\hat x_1)$.
\item $\forall v \in E_i(\hat x) ^*$, $\lim {1 \over n} \log \abs{d_0 f_x^n(v)} = \Lambda_i$, and $\forall v \in E_{u+1}(\hat x) ^*$, $\lim {1 \over n} \log \abs{d_0 f_x^n(v)} \leq 0$.
\end{enumerate}
 There exist $C : \orb \to GL_k(\C)$ and an $\epsi$-tempered function $\varphi_\epsi$ such that:
\begin{enumerate}
\item[4.] $C_{\hat x}$ sends $\oplus_{i=1}^{u} \C^{m_i} \oplus \C^p$  to $\oplus_{i=1}^u E_i(\hat x) \oplus E_{u+1}(\hat x)$ and $\abs v \leq \abs {C_{\hat x}(v)} \leq \varphi_\epsi(\hat x) \abs{v}$.
\item[5.]  $D_{\hat x} := C_{\hat f (\hat x)} \circ d_0 f_x \circ C_{\hat x}^{-1}$ is a block diagonal map $(D_{\hat x}^1 , \ldots , D_{\hat x}^{u+1} )$.
\item[6.]  $\forall v \in \C^{m_i}$, $e^{\Lambda_i - \epsi} \abs{v} \leq \abs {D_{\hat x}^i (v)} \leq  e^{\Lambda_i + \epsi} \abs{v}$ and $\forall v \in \C^p$, $\abs{D_{\hat x}^{u+1} (v)} \leq
  e^{ \epsi} \abs{v}$.
\end{enumerate}
\end{thm}

We set $\zeta_{\hat x} := \tau_{x_0} \circ C_{\hat
  x}^{-1}$ and $g_{\hat x} := \zeta_{\hat x_1}^{-1} \circ f \circ \zeta_{\hat
  x}$. Observe that $\zeta_{\hat x} : \D^k(\rho_0) \to
  \Pj^k$ and $g_{\hat x}(0) = 0$, $d_0 g_{\hat x} =
D_{\hat x}$. The following lemma holds up to multiply $\varphi_\epsi$
  and divide $\rho_0$ by a multiple of $1 + \abs {d f}_\infty
  + \abs {d^2 f}_\infty$ depending on the derivatives of the $\tau_x^{\pm 1}$'s. 

\begin{lem}\label{poli}
For every $r \leq \rho_0$, we have:
\begin{enumerate}
\item $\forall (u,v) \in \D^k(r)$, $\varphi_\epsi ^{-1}(\hat x) \abs{u-v} \leq d(\zeta_{\hat x}(u),\zeta_{\hat x}(v)) \leq  2 \abs{u-v}$.
\item The map $g_{\hat x}$ is well defined from $\D^k(r)$ to $\D^k(r \varphi_\epsi(\hat x_1))$.
\item $\forall w \in \D^k(r)$, $\abs{ d_w g_{\hat x} - d_0 g_{\hat x} } \leq r \varphi_\epsi (\hat x_1) $.
\end{enumerate}
\end{lem}

The proof of that lemma is left to the reader. We set for the sequel $\varphi_0 \geq 1$  such that $\widehat \Omega := \{ \varphi_\epsi \leq
 \varphi_0 \}$ satisfies $\hat m (\widehat \Omega) \geq 9/10$. Now our aim is to introduce the entropy of $m$. Let $d_n$ be the  distance on $\Pj^k$ defined by $d_n(x,y) := \max_{0 \leq j \leq n-1} d(f^j(x),f^j(y))$ and let $B_n(x,r) := \{ y \in \Pj^k \, , \,  d_n(x,y)  <  r  \}$. Brin-Katok's theorem \cite{BK}
 asserts that for
$m$-a.e. $x\in \Pj^k$, we have: 
\[ \sup_{r > 0} \,  \liminf_{n \to +\infty} \,  - {1 \over n} \log m (B_n(x,r)) = h_m . \] 
We shall only need: there exist $r_\epsi(x) > 0$ and $n_\epsi(x) \geq 0$ satisfying
\begin{equation}\label{entbk}
\forall  r \leq  r_\epsi(x) \, , \, \forall  n \geq n_\epsi(x) \, , \,
m(B_n(x,r)) \leq e^{-n(h_m - \epsi)} .
\end{equation}
We set $r_0 > 0$ and $n_0 \geq 1$ such that $\BB := \{ r_\epsi
\geq r_0 \, , \, n_\epsi  \leq n_0  \} $ satisfies
$m(\BB) > 9/10$.  For every $n \geq
n_0$, we fix a maximal $r_0$-separated subset $\EE_n$ (for the distance $d_n$)
in $\pi(\widehat \Omega) \cap
\BB$. We have $\cdb \EE_n \geq e^{n(h_m-2\epsi)}$ by using (\ref{entbk}).

%%%%%%%%%%%%%%%%%%%
\subsection{Proof of theorem D}
%%%%%%%%%%%%%%%%%%%

The proof relies on the following proposition (we set $q_{u+1} := 0$). We recall that it suffices to establish theorem D for $2 \leq j \leq u+1$. 

\begin{prop}\label{danm}
Let $n \geq n_0$, $x \in \EE_n$ and $\ga_0 \leq 1$. For every $2 \leq j \leq u+1$, there exists a neighbourhood $U_x^j$ of the origin in  $\D^{p+q_j}$ and a mapping $\Psi^j_x : U_x^j \to \Pj^k$  (depending on $n$) which satisfies the following properties:

\begin{enumerate}

\item $\Psi^j_x(0) = x$ and $\Lip \Psi^j_x \leq \ga_0$.

\item $\Diam f^i
  (\Psi^j_x) \leq e^{-n\epsi}$ for every $0 \leq i \leq n-1$.

\item  $\Vol \ \Psi^j_x \geq e^{- n(2m_{j} \Lambda_{j} + \ldots +
  2m_u \Lambda_u)} e^{-8 k n \epsi}$ for $2 \leq j \leq u$, and $\Vol \ \Psi^{u+1}_x \geq e^{-8 k n \epsi}$.

\end{enumerate}

\end{prop}

Let us assume  proposition \ref{danm} and complete the proof of theorem D.  The following arguments were employed in \cite{dT}, we give them for
  reader's convenience. Let $2 \leq j \leq u+1$ and $\Psi^j := \cup_{x
  \in \EE_n}  \Psi_x^j$. Taking affine
charts, we can work on
$\C^k$ and assume that the $\Psi^j_x$'s are graphs above $\PP = 
\C^{p+q_j}$. Let $\si : \C^k \to \PP$ be the orthogonal projection and $\si_a := \si^{-1}\{ a\}$.  Let $\om$
  be the Fubini-Study form on $\Pj^k$ and endow $ \Pj^k \times \ldots \times \Pj^k$ ($n$ times) with the metric $\om_n := \sum_{i=1}^n
p_i^* \om$, where the $p_i$'s denote the projections to the factors. Denoting $\Ga_n(a) := \{ (
z , f(z) , \ldots , f^{n-1}(z)  ) \, , \, z \in \si_a \}$, we have $\Vol \, \Ga_n(a)  = \int_{\Ga_n(a)} \om_n^{k-p-q_j}$.
The crucial observation is that $\Vol \ \Ga_n(a) \geq \cdb \Psi^j \cap
\si_a$ for every $a \in \AA$. This is a consequence of  Lelong's
inequality (the holomorphic context is crucial here) and the fact that $\Psi^j
 \cap \si_a$ is $r_0/2$-separated for the distance $d_n$ (that follows from the fact that $\EE_n$ is $r_0$-separated and from $\Diam f^i
  (\Psi^j_x) \leq e^{-n\epsi}$ provided by proposition \ref{danm}(2)). That observation implies after an integration over $a \in \PP$:
 
  \begin{equation}\label{volpro}
  \int_{a \in \PP} \Vol \, \Ga_n(a) \, da \geq \int_{a \in \PP} \cdb \Psi^j \cap \si_a \, da = \Vol \,
 \si(\Psi^j).
 \end{equation}
Using proposition \ref{danm}(3), we deduce $(\star)$: $\Vol \, \si(\Psi^j)  \geq   e^{n (h_m - 2\epsi)}  e^{- n (2m_{j} \Lambda_{j} + \ldots +
  2m_u \Lambda_u)} e^{-8kn\epsi}$ (we deal with $2 \leq j \leq u$). We now focus on the upper bound. Let $[\si_a]$ be the current of
integration on $\si_a \simeq \C^{k-p-q_j}$ and $\Omega := \int_{a \in \PP}  [\si_a] \, da$. We have by definition of $\om_n$:
\[  \int_{a \in \PP}  \Vol \, \Ga_n(a) \, da  =  \sum \ \int_{\Pj^k}
\Omega \wedge {f^{n_1}}^* \om \wedge \ldots \wedge
   {f^{n_{k-p-q_j}}}^* \om  , \]
where the sum runs over $0 \leq n_1, \ldots, n_{k-p-q_j} \leq n-1$. Since $\Omega \leq \om^{p+q_j}$, we get:
\[   \int_{a \in \PP}  \Vol \, \Ga_n(a) \, da  \leq \sum \ \int_{\Pj^k}
\om^{p+q_j} \wedge {f^{n_1}}^* \om \wedge \ldots \wedge
   {f^{n_{k-p-q_j}}}^* \om  . \]
Now the inner integral is equal to $d^{n_1 + \ldots + n_{k-p-q_j}}$
   since ${f^n}^*\om$ is cohomologous to $d^n \om$ (see \cite{DS2}, section 1.2). We deduce $(\star \star)$: $\int_{a
   \in \PP}  \Vol \, \Ga_n(a) \, da \leq n^{k-p-q_j} d^{n(k-p-q_j)}$. Combining (\ref{volpro}), $(\star)$ and $(\star \star)$, we finally obtain:
\[   \log d^{ k-p-q_j } +  n^{-1} \log n^{k-p-q_j}  \geq   h_m   -  (2m_j\Lambda_j + \ldots +
  2m_u \Lambda_u) -(8k+2) \epsi .  \]
Theorem D then follows by taking limits and observing $k-p-q_j = m_1 + \ldots + m_{j-1}$.

%%%%%%%%%%%%%%%%%%%
\subsection{Proof of proposition \ref{danm}}
%%%%%%%%%%%%%%%%%%%

We shall work in the charts $\zeta_{\hat x} : \D^k(\rho_0) \to \Pj^k$ (see subsection
\ref{ospe}), that is with the local maps $g_{\hat x} = \zeta_{\hat x_1}^{-1} \circ f \circ \zeta_{\hat
  x}$. Let us prove the following assertions  for every $\hat x \in \widehat
\Omega$. In what follows $i$ ranges $\{  0, \ldots , n-1\}$ and the $\psi_i^j$'s depend on $n$.

\begin{enumerate}

\item[($\al$)] Let $2 \leq j \leq u$. There exist $\psi_i^j :
U_i \subset \D^{p+q_j}(e^{-3n\epsi})\to \C^{k-p-q_j}$ such that: 

- $\psi_i^j(0)=0$ and $\Lip \psi_i^j \leq \ga_0$,

- $g_{\hat x_i} (\graphe
\psi_i^j) \subset \graphe \psi_{i+1}^j$,

- $\Vol \ \psi_0^j \, \geq  \, e^{- n(2m_{j} \Lambda_{j} + \ldots +
  2m_u \Lambda_u)} e^{-8kn\epsi}$.

\item[($\beta$)] The same properties hold for $j=u+1$, with $\psi_i^{u+1} : \D^p(e^{-3n \epsi - (n-i)\epsi})
  \to \C^{k-p}$ and the last item being replaced by  $\Vol \ \psi_0^{u+1}  \, \geq  \, e^{-8kn \epsi}$.
 
\end{enumerate} 

Proposition \ref{danm} follows from ($\al$) and ($\beta$) by setting $\Psi^j_x :=
 \zeta_{\hat x} \circ \psi_0^j$ and using $\varphi_0^{-1} \abs{u-v} \leq d(\zeta_{\hat x}(u),\zeta_{\hat x}(v)) \leq  2 \abs{u-v}$ (see lemma \ref{poli}). The points
 1 and 3 of that proposition are clear up to multiplicative constants. The point 2 is a consequence of $f^i(\Psi^j_x) =  \zeta_{\hat
 x_i} (\graphe \psi_i^j)$, $U_i \subset \D^{p+q_j}(e^{-3n\epsi})$ and $\Lip \psi_i^j  \leq 1$. \\
 
 The proofs of
 ($\al$) and ($\beta$) rely onbackward graph
 transforms for possible non injective maps $g$ (see theorem \ref{bgt1} below, the partial
 derivative $B$ can be zero). We will successively apply it for
 $g_{\hat x_i}$ from $i = n-1$ to $i=0$. Observe  that theorem \ref{newhouse}(6) ensures to satisfy conditions (a)-(b) for $A^l_{\hat x} := (D_{\hat x}^1,\ldots, D_{\hat
 x}^l)$, $B^l_{\hat x} := (D_{\hat x}^{l+1},\ldots, D_{\hat
 x}^{u+1})$, $1 \leq l \leq u$, and conditions (c)-(d) for $(A^{u}_{\hat x}$, $B^{u}_{\hat x})$. Moreover, given a small $\delta < \epsi$ and setting $R_0 = e^{-3n\epsi}$, $R_1 = e^{-3n\epsi}\varphi_\epsi(\hat x_{i+1})$,  lemma 
 \ref{poli}(3) ensures that $g_{\hat x_i} : \D^k(R_0) \to \D^k(R_1)$ satisfies $\abs{ d_w g_{\hat x_i} - d_0 g_{\hat x_i} }
 \leq  \delta$: we indeed have $\abs{ d_w g_{\hat x_i} - d_0 g_{\hat x_i} }
 \leq e^{-3n\epsi} \varphi_0 e^{(i+1)\epsi}  \leq \varphi_0
 e^{-2n_0\epsi}$ for every  $0 \leq i \leq n-1$ and $n \geq n_0$. The backward graph transform is stated as follow, the proof is postponed to subsection \ref{pbgt}.
 
\begin{thm}\label{bgt1} Let $(k_1 , k_2)$ be positive integers such that $k = k_1
  + k_2$, and $A : \C^{k_1} \to \C^{k_1}$, $B : \C^{k_2} \to \C^{k_2}$ be linear maps. Assume that $A$ is invertible
  and that $\abs B < \abs{ {A}^{-1} }^{-1}$. We denote $\ga := 1 -
  \abs B \abs {A^{-1}} \in ]0,1]$. Let $0 \leq \ga_0 \leq
  1$ and $\delta < \epsi$ be such that  
\begin{enumerate}
\item[(a)] $\ga_0(1-\ga) + 2 \delta(1+\ga_0) \abs {A^{-1}} \leq 1$,
\item[(b)] $\left( \ga_0 \abs B + \delta (1+\ga_0) \right) \left( \abs
  {A^{-1}}^{-1} - \delta (1+\ga_0) \right)^{-1} \leq \ga_0$.
\end{enumerate}
Let  $g : \D^k(R_0) \to \D^k(R_1)$ be an holomorphic mapping such that $R_0 \leq R_1$, $g(0)=0$, $d_0 g = (A,B)$ and $\abs{d_w g - d_0 g} \leq \delta$ on $\D^k(R_0)$. 

\begin{enumerate}

\item  If $\phi : V \subset \D^{k_2}(R_1) \to \C^{k_1}$ satisfies
  $\phi(0) = 0$ and $\Lip \phi \leq \ga_0$, then there exists $\psi : U \subset
\D^{k_2}(R_0) \to \C^{k_1}$ such that
\begin{equation}\label{incc}  
\Lip \psi \leq \ga_0 \ \textrm{  and } \   g  (\graphe \psi)  \subset \graphe \phi.
 \end{equation}

\item Assume (c): $\abs B + 2\delta \leq
  e^{\epsi}$. If $\phi : \D^{k_2}(R)  \to \C^{k_1}$ satisfies
  $\phi(0) = 0$ and $\Lip \phi \leq \ga_0$ for some $R \leq R_0$, then there exists $\psi : \D^{k_2}(R e^{-\epsi}) \to \C^{k_1}$ satisfying (\ref{incc}).

\item  Assume (d): \[ (\abs B + 2\delta) e^{-\epsi}  + \delta \leq 1  \ \ \textrm{and} \ \ \delta (1 + \ga_0)  \leq   \min \{   ( \abs {A^{-1}}^{-1} - \ga_0 \abs B)  / 2,  \abs {A^{-1}}^{-1} - 1 \} . \]
If $\phi : \D^{k_2}(R) \to \C^{k_1}$ satisfies  $\abs {\phi(0)} \leq
R$ and $\Lip \phi \leq \ga_0$ for some $R \leq R_0/2$, then there exists $\psi : \D^{k_2}(Re ^{-\epsi}) \to \C^{k_1}$ satisfying $\abs {\psi(0)} \leq R$ and (\ref{incc}).
\end{enumerate}
\end{thm}

In order to prove $(\beta)$, we apply theorem \ref{bgt1}(2) starting with $\phi := \psi_n : \D^p(e^{-3n\epsi}) \to
\C^{k-p}$ equal to zero and splitting along the orbit according to
$(A^{u}_{\hat x_{i}}, B^{u}_{\hat x_{i}})$. We obtain mappings $\psi_i :
\D^p(e^{-3n\epsi -(n-i)\epsi}) \to \C^{k-p}$ such that $\Lip \psi_i \leq \ga_0$
and $g_{\hat x_i} (\graphe \psi_i)  \subset \graphe \psi_{i+1}$. We
obviously have in that case $\Vol \, \psi_0  \geq e^{-8kn \epsi}$ since $\psi_0$ is defined on $\D^p(e^{-4n\epsi})$.\\

For ($\alpha$), we apply theorem  \ref{bgt1}(1) starting with $\psi_n : \D^{p+q_j}(e^{-3n\epsi}) \to
\C^{k - p - q_j}$ equal to zero and splitting according to $(A^{j-1}_{\hat x_i},B^{j-1}_{\hat x_i})$. We obtain mappings $\psi_i :
U_i \subset \D^{p+q_j}(e^{-3n\epsi}) \to \C^{k-p-q_j}$ such that $\Lip \psi_i \leq \ga_0$
and $g_{\hat x_i} (\graphe \psi_i)  \subset \graphe \psi_{i+1}$. 

It remains to estimate the volume of $\psi_0$. We use for that purpose the slicing
argument of \cite{dT}, the idea is to control the size of the $U_i$'s. Let us write $\psi_n = \cup_{a \in \D^{q_j}(e^{-3n\epsi})} \,  \psi_{a,n}$, where $\psi_{a,n} : \D^p(e^{-3n\epsi}) \to
\C^{k - p}$ is equal to $\{ 0 \}^{k-p-q_j} \times \{ a
\}$. Now apply theorem \ref{bgt1}(3) starting with $\phi =
\psi_{a,n}$ and $R = e^{-3n\epsi}$, $R_0 = e^{-2n\epsi}$. We obtain mappings $\psi_{a,i} :
\D^{p}(e^{-3n \epsi - (n-i) \epsi}) \to \C^{k-p}$ satisfying $\psi_i
= \cup_a \psi_{a,i}$. In
particular, $\psi_0$ is foliated by graphs above
$\D^p(e^{-4n\epsi})$. 

Let $\psi_0^b :=
\psi_0 \cap ( \D^{k-p}
\times \{ b \})$ for  $b \in \D^p(e^{-4n\epsi})$, and let $\psi_i^b := g_{\hat x_{i-1}} \circ \ldots \circ g_{\hat
  x} (\psi_0^b)$. Since $\psi_n^b$ intersects every $\psi_{a,n}$, we have $\Vol \, \psi_n^b \geq e^{-6 q_j n
  \epsi}$. The next step is to show: 
 \begin{equation}\label{abn} \forall b \in \D^p(e^{-4n\epsi}) \ ,  \ \Vol \, \psi_0^b \ e^{n( 2m_{j} \Lambda_{j} +
    \ldots + 2m_u \Lambda_u)  + 2q_j n  \epsi} \geq \Vol \, \psi_n^b.
    \end{equation}
This can be done inductively by checking $\abs { \wedge^{q_j}   {d_z
    g_{\hat x_i}}} \leq e^{2m_{j} \Lambda_{j} +
    \ldots + 2m_u \Lambda_u + 2q_j \epsi}$ on $\psi_i^b$, we refer to \cite{dT} for the details. We observe that the estimate (\ref{abn}) holds even if multiplicities occur for the
    jacobians. Finally, the combination of the coarea formula 
   $\Vol \, \psi_0 = \int_{\D^p(e^{-4n\epsi})} \, \Vol \, \psi_0^b \ db$, (\ref{abn}) and $\Vol \, \psi_n^b \geq e^{-6 q_j n
  \epsi}$ yields item ($\alpha$).

%%%%%%%%%%%%%%%%%%%%%%%%%%%%%%%%%%%%%%%%%%%%%%%%%%%%
\subsection{Proof of the backward graph transform}\label{pbgt}
%%%%%%%%%%%%%%%%%%%%%%%%%%%%%%%%%%%%%%%%%%%%%%%%%%%%

We  follow the method of Hirsch-Pugh-Shub (see \cite{HPS}, section 5). We denote by $(x,y)$ the elements of $\C^{k_1} \times  \C^{k_2}$ and $g = (g_1,g_2)$. Let us prove the point 1.  Let $\phi : V \subset \D^{k_2}(R_1) \to \C^{k_1}$ satisfying $\phi(0) = 0$ and $\Lip \phi \leq \ga_0$. For any $y \in \D^{k_2}(R_0)$, let $L(y) := \{ (x,y) \, , \, \abs x
\leq \abs {y} \} \subset \D^k(R_0)$ and 
\[ U:= \{ y \in \D^{k_2}(R_0)
\, , \,  g_2(L(y)) \subset V \}. \] 
We fix $y \in U$ and find $(x,y) \in L(y)$ such that
$\phi (g_2(x,y)) = g_1(x,y)$. This is equivalent to find $(x,y) \in L (y)$ such that:
\[ \Lambda_{y}(x) := A^{-1} \left[  \phi(g_2(x,y)) - (g_1(x,y) - Ax)  \right] = x .    \] 
First let us verify that $\Lambda_{y} : \bar \D^{k_1}(\abs{y}) \to
\bar \D^{k_1}(\abs{y})$. Using $g(0) = 0$, $\abs{d_w g - d_0 g} \leq \delta$ and $d_0 g = (A,B)$ on one hand, and $\phi(0)=0$, $\Lip \phi \leq \ga_0$ on the other hand, we obtain
\begin{displaymath}
\begin{array}{rcl}
 \abs {\Lambda_{y}(x) }  & \leq  &  \abs{A^{-1}} \big[ \abs { \phi
    (g_2(x,y)) - \phi(0) } + \abs{ g_1(x,y) - Ax
    }   \big]       \\
                                                    & \leq & \abs{A^{-1}} \big[ \ga_0 (\delta
    \abs x +  (\abs B + \delta)\abs{y} ) + \delta \abs x + \delta \abs{y}  \big]  .
\end{array}                                      
\end{displaymath}
Using $\abs x \leq \abs{y}$, we have $\abs {\Lambda_{y}(x) } \leq
\abs{A^{-1}} \big[ \ga_0 \abs B +   2\delta (1+\ga_0)  \big ]
\abs{y}$. From $\abs {A^{-1}} \abs B = 1 -\ga$ and condition  (a), we obtain  $\Lambda_y : \bar \D^{k_1}(\abs{y}) \to \bar \D^{k_1}(\abs{y})$. Let us now prove that $\Lambda_{y}$ is contracting. The difference $\abs {\Lambda_{y}(x) - \Lambda_{y}(x') }$ is less than or equal to 
\begin{displaymath}
\begin{array}{rcl}
 & \abs{A^{-1}}
\Big[ \abs { \phi(g_2(x,y)) -  \phi(g_2(x',y)) }  +  \abs{
    (g_1(x,y) - Ax) - (g_1(x',y) - Ax' )}  \Big] &  \\
         & \leq \abs{A^{-1}}
\Large[ \ga_0 \delta \abs{x-x'}  + \delta \abs{x-x'} \Large] . &  
\end{array}                                      
\end{displaymath}
The contraction property then follows from (a) (it implies $\delta (1+\ga_0) \abs {A^{-1}} \leq 1/2$). For every $y \in U$,
we denote by $\psi(y)$ the unique fixed point of $\Lambda_{y}$. Let us show that $\Lip \psi \leq \ga_0$. For
every $(y,y') \in U$, we have $\abs{\psi(y) - \psi(y')} = \abs{\Lambda_y (\psi(y) )
- \Lambda_{y'} (\psi(y') )}$. Hence $\abs{\psi(y) - \psi(y')} \leq  \abs{A^{-1}} (I + J)$,  
where 
\[ I :=  \abs { \phi
    (g_2(\psi(y),y)) - \phi
    (g_2(\psi(y'),y')) } \leq \ga_0 \left[ \delta \abs{\psi(y) - \psi(y')}
    + (\abs B + \delta)\abs{y-y'} \right], \] 
\[  J := \abs{ \left[ g_1(\psi(y),y) - A(\psi(y)) \right] - \left[ g_1(\psi(y'),y') -
    A(\psi(y')) \right] } \leq \delta \abs{\psi(y) - \psi(y')} +
    \delta \abs{y-y'}   .  \]  
We deduce:
\[ \abs{\psi(y) - \psi(y')} \leq { \ga_0 \abs B +
    \delta (1+\ga_0)  \over  \abs
  {A^{-1}}^{-1} - \delta (1+\ga_0)  }   \abs{y-y'}, \]
which is less than $\ga_0 \abs{y-y'}$ by using (b). That proves the
point 1. For the  point 2, observe that for every $y \in \D^k(R e^{-\epsi})$ and $(x,y) \in L(y)$ (use (c)):  
\[ \abs{g_2(x,y)} \leq \delta \abs x + (\abs B + \delta)\abs {y}
\leq   (\abs B + 2 \delta) \abs{y} < R . \]
That implies $\D^{k_2}(R e^{-\epsi}) \subset U := \{ g_2(L(y)) \subset \D^{k_2}(R) \}$. The point 2 then follows by repeating the same arguments than for the point 1. Now the point 3. Let $y \in \D^{k_2}(R e^{-\epsi})$ and $L_R(y) := \{ (x,y) \, , \, \abs x
\leq \abs {y} + R \}$. Observe that $L_R(y) \subset \D^k(R_0)$, because $R e^{-\epsi} + R \leq 2R \leq R_0$. As before we look for $(x,y) \in L_R(y)$ such that
$\phi (g_2(x,y)) = g_1(x,y)$. Observe that for every $(x,y) \in L_R(y)$ (use (d) for the last estimate):
\[ \abs{g_2(x,y)} \leq \delta \abs x + (\abs B + \delta)\abs {y}
\leq (\abs B + 2 \delta) \abs{y} + \delta R <  [ (\abs B + 2 \delta)  e^{-\epsi} + \delta] R \leq R . \]
That yields $\D^{k_2}(R e^{-\epsi}) \subset U$. Let us show that $\Lambda_{y} : \bar \D^{k_1}(\abs{y}+R) \to
\bar \D^{k_1}(\abs{y}+R)$. Recalling that  $\Lambda_{y}(x) = A^{-1} \left[
   \phi (g_2(x,y))  - ( g_1(x,y) - Ax ) \right]$, we have:
\begin{displaymath}
\begin{array}{rcl}
 \abs {\Lambda_{y}(x) }  & \leq  &  \abs{A^{-1}} \big[ \abs { \phi
    (g_2(x,y)) - \phi(0) } + \abs{\phi(0)} + \abs{ g_1(x,y) - Ax
    }   \big]          \\
                                                    & \leq & \abs{A^{-1}} \big[ \ga_0 \Large( \delta
    \abs x + ( \abs B + \delta ) \abs{y} \Large)  + R + ( \delta \abs x
    + \delta \abs{y})  \big]  .
\end{array}                                      
\end{displaymath}
We deduce using $\abs x \leq \abs{y} + R$:
\[ \abs {\Lambda_{y}(x) } \leq  \abs{A^{-1}} [  \ga_0 \abs B + 2 \delta (1+\ga_0)  ] \abs{y}      +  \abs{A^{-1}}   ( 1 + \delta (1+\ga_0)) R.   \]
It follows from (d) that $\abs {\Lambda_{y}(x) } \leq
\abs{y} + R$. The same arguments as before give the contraction property for $\Lambda_{y}$ and
the Lipschitz property of $\psi$ ($\phi(0)$ is not
involved here). Observe finally that $\abs{\psi(0)} \leq R$, because 
$\Lambda_0 : \bar \D^{k_1}(R) \to \bar \D^{k_1}(R)$.

\vspace{.6 cm}

{\footnotesize C. Dupont}\\
{\footnotesize Universit\'e Paris XI}\\
{\footnotesize CNRS UMR 8628}\\
{\footnotesize Math\'ematiques, B\^atiment 425}\\
{\footnotesize F-91405 Orsay Cedex, France}\\
{\footnotesize christophe.dupont@math.u-psud.fr}

\end{document}